\documentclass{amsart} \textwidth=14.5cm \oddsidemargin=1cm
\usepackage{amsmath}
\usepackage{amsxtra}
\usepackage{amscd}
\usepackage{amsthm}
\usepackage{amsfonts}
\usepackage{amssymb}
\usepackage{eucal}
\usepackage{epigraph}

\newtheorem{Thm}{Theorem}
\newtheorem{Cor}{Corollary}
\newtheorem{Lem}{Lemma}
\newtheorem{Prop}{Proposition}
\newtheorem{Sublem}{Sublemma}

\newtheorem{Claim}{Claim}
\theoremstyle{remark}
\newtheorem{Rem}{Remark}

\newcommand{\Ql}{{\overline {{\Bbb Q}_l} }}

\newcommand{\cal}{\mathcal}
\newcommand{\Gr}{{{\cal G}{\frak r} }}
\newcommand{\Fl}{{{\cal F}\ell}}

\newcommand{\dfn}{{      \overset{\rm {def}} = }} 
\newcommand{\bu}{\bullet}
\newcommand{\fr}{fr}
\newcommand{\To}{\longrightarrow}
\newcommand{\iso}{{\widetilde \longrightarrow}}

\newcommand{\imbed}{\hookrightarrow}

\newcommand{\<}{\langle}
\renewcommand{\>}{\rangle}
\newcommand{\de}{\partial}
\def\square{\hbox{\vrule\vbox{\hrule\phantom{o}\hrule}\vrule}}
\newcommand{\epf}{\square}
\newcommand{\unV}{{\underline V}}

\newcommand{\jbar}{{\bar j}}

\renewcommand{\P}{{\cal P}}
\newcommand{\PI}{{\cal P}_I}
\newcommand{\fP}{{^f{\cal P}}}

\newcommand{\Po}{{{\cal P}^0}}
\newcommand{\PoI}{{{\cal P}^0_I}}
\newcommand{\PIW}{\cal P_{{\cal {IW}}}}
\newcommand{\FIW}{F_{{\cal {IW}}}}
\newcommand{\DIW}{D_{{\cal {IW}}}}

\newcommand{\io}{{\mathfrak i}}
\newcommand{\Si}{{\mathfrak S}}

\newcommand{\cS}{{\cal S}}

\renewcommand{\H}{{\cal H}}
\newcommand{\fPhi}{{^f{\Phi}}}

\newcommand{\Ftil}{\tilde{F}}

\newcommand{\N}{{\cal N}}
\newcommand{\Ntil}{{\tilde{\cal N}}}

\newcommand{\Nt}{{\tilde{\cal N}}}

\newcommand{\wtil}{{\tilde w}}
\newcommand{\lambdatil}{{\tilde \lambda}}

\newcommand{\suml}{\sum\limits}
\newcommand{\oplusl}{\bigoplus\limits}
\newcommand{\cupl}{\bigcup\limits}
\newcommand{\fW}{{^f W}}

\newcommand{\fWf}{{^fW^f}}

\renewcommand{\b}{{\frak b}}

\newcommand{\bo}{B}
\newcommand{\bobo}{b}
\newcommand{\bobobo}{{\bf b}}
\newcommand{\bFl}{{\bf b}}
\newcommand{\bq}{{\bf q}}

\newcommand{\M}{{\cal M}}
\newcommand{\mon}{{\mathfrak M}}

\newcommand{\Lg}{{ \frak g\check{\ }}}
\newcommand{\Ln}{{ \frak n\check{\ }}}
\newcommand{\LG}{{ G\check{\ }}}
\newcommand{\LB}{{B\check{\ }}}
\newcommand{\LP}{{P\check{\ }}}
\newcommand{\LU}{{U\check{\ }}}
\newcommand{\LT}{{T\check{\ }}}

\newcommand{\LGU}{{ G\check{\ }/U\check{\ }}}
\newcommand{\LGUbar}{\overline{ G\check{\ }/U\check{\ }}}
\newcommand{\BF}{{\bf B_F}}
\newcommand{\NFm}{{\bf N^-_F}}
\newcommand{\BFm}{{\bf B^-_F}}
\newcommand{\NF}{{\bf N_F}}

\newcommand{\Nthat}{\hat{\tilde{\cal N}}}

\renewcommand{\L}{{\cal L}}

\renewcommand{\O}{{\cal O}}
\newcommand{\Ohat}{{\hat\O}}

\newcommand{\F}{{\cal F}}
\newcommand{\G}{{\cal G}}
\newcommand{\A}{{\cal A}}
\newcommand{\B}{{\cal B}}
\newcommand{\C}{{\cal C}}
\newcommand{\E}{{\cal E}}
\newcommand{\cT}{{\cal T}}
\newcommand{\grA}{gr{\cal A}}

\newcommand{\cons}{{\overline{\underline{{\Bbb Q}_l}}}}
\newcommand{\Z}{{\cal Z}}

\newcommand{\Zet}{{\Bbb Z}}
\newcommand{\Ce}{{\Bbb C}}

\newcommand{\Fbar}{{\bar F}}
\newcommand{\Lotimes}{\overset{\rm L}{\otimes}}
\newcommand{\unO}{{\underline{\cal O}}}

\newcommand{\bI}{{\bf I}}

\newcommand{\GO}{{\bf G_O}}
\newcommand{\K}{\GO}
\newcommand{\GK}{{\bf G_F}}
\newcommand{\GF}{{\bf G_F}}
\newcommand{\BO}{{\bf B_O}}
\newcommand{\BOm}{{\bf B^-_O}}

\newcommand{\Db}{D^b}

\newcommand{\bbA}{{\mathbb A}}
\newcommand{\bbH}{{\mathbb H}}
\newcommand{\bbK}{{\mathbb K}}
\newcommand{\bbM}{{\mathbb M}}

\newcommand{\Aone}{{\mathbb A}^1}
\newcommand{\Ga}{{\mathbb G}_a}
\newcommand{\Gm}{{\mathbb G}_m}
\newcommand{\AS}{{\mathbb {AS}}}
\newcommand{\Fq}{{\mathbb F}_q}

\newcommand{\Fpbar}{{\bar{\mathbb F}_p}}
\newcommand{\Fpn}{{\mathbb F}_{p^n}}
\newcommand{\lT}{\frac{\ell(w)}{2}}

\newcommand{\AutO}{{\bf Aut}(O)}

\newcommand{\Sa}{S}

\newcommand{\Stalk}{Stalk}
\newcommand{\Costalk}{Costalk}

\newcommand{\la}{\lambda}
\newcommand{\La}{\Lambda}

\renewcommand{\proof}{{\it Proof }}
\newcommand{\proofpt}{{\it Proof. }}

\author{Sergey Arkhipov and Roman Bezrukavnikov}
\title[Affine flags and the dual group]
{Perverse sheaves on affine flags and Langlands dual group}

\begin{document}
\maketitle

\medskip 

\epigraph{Though he goeth on his way weeping that beareth the measure of seed, he 
shall  come home with joy, bearing his sheaves.}
{Psalm 126 "Shir ha-maalot" (Song of Ascendance)}

\tableofcontents

\begin{abstract} This is the first in a series of 
papers devoted to describing the category of sheaves on the affine
flag manifold of a simple algebraic group in terms of the Langlands
dual group. In the present paper we provide such a description for
categories which are geometric counterparts of a maximal
commutative subalgebra in the Iwahori Hecke algebra $\bbH$; of the
anti-spherical module for $\bbH$; and of
 the space of
Iwahori-invariant Whittaker functions. As a byproduct we obtain
some new properties of central sheaves introduced in \cite{KGB}.
\end{abstract}

\thanks{{\bf Acknowledgements.}
 This project was conceived during the IAS special year
in Representation Theory (1998/99) led by G.~Lusztig,
as a result of conversations
with D.~Gaitsgory, M.~Finkelberg and I.~Mirkovi\' c. The outcome was
strongly influenced by conversations with A.~Beilinson and V.~Drinfeld.
The stimulating interest of A.~Braverman, D.~Kazhdan,
G.~Lusztig and V.~Ostrik was crucial
for keeping the project alive. 
We are very grateful to all  these people. We  thank
 I.~Mirkovi\' c
and D.~Gaitsgory for the permission to use their unpublished
results; and 
M.~Finkelberg and D.~Gaitsgory for taking the trouble to read
the text and point out various lapses in the exposition. 
S.A. was partly supported by an NSERC grant, and R.B.
 by an NSF grant
and the Clay Institute while working on this paper.}

\begin{section}{Introduction}

\subsection{Basic notations and motivation}
 Let $k$ be a field of characteristic $p>0$; we will soon set $k=\Fpbar$,
but in the introduction we also allow $k$ to be finite;
 $G$    % \supset B=T\cdot N$
 be a split simple linear algebraic group over $k$
 with a Borel subgroup $B$, $k((t))=F\supset
O=k[[t]]$
be a local functional field  and its ring of integers.
 Let $G(O)\supset I$ be respectively a maximal compact
subgroup, and an Iwahori. %, and its unipotent radical.

 There exist
canonically defined group schemes $\K$, $\bI$ over $k$ (of
infinite type) such that  $\K(k)=G(O)$, $\bI(k)=I$; and an
ind-group scheme $\GF$ with $\GF(k)=G(F)$. We also have the
homogeneous
ind-varieties:   the affine flag variety $\Fl$ and
the affine Grassmanian $\Gr$, see e.g. \cite{KGB}, Appendix, \S A.5.
 Thus $\Fl$, $\Gr$ are direct limits of
projective varieties with transition maps being closed embeddings,
and $\Fl(k)=G(F)/I$, $\Gr(k)=G(F)/G(O)$.

%Let $\pi:\Fl\to\Gr$ be the projection.
 %If $k$ is algebraically closed,
  Let $D=D(\Fl)$, $D(\Gr)$ be the constructible 
 derived category of $l$-adic sheaves ($l\ne
char(k)$; see \cite{Weil2}, 1.1.2; %check 
\cite{BBD} 2.2.14--2.2.18; and also \cite{KGB}, \S A.2) 
 on
 $\Fl$, $\Gr$ respectively, and
$D_I=D_I(\Fl)$, $D_{G(O)}=D_{G(O)}(\Gr)$ be the equivariant derived
categories (cf. \cite{BL}).  Let
 $\P\subset D$, $\P(\Gr)\subset D(\Gr)$,
$\P_I\subset D_I$, $\P_{G(O)}\subset D_{G(O)}$ be the subcategories of
perverse sheaves.

 By $*$ we
denote the convolution; thus $*$ provides $D_{G(O)}(\Gr)$, $D_I(\Fl)$,
with a monoidal structure,
 and defines a "right"
action of $D_I(\Fl)$ 
on $D(\Fl)$.

 Let $\LG$ be  the
Langlands dual group
over the field $\Ql$, and $Rep(\LG)$ be its category of representations.

Recall that according to a result of Lusztig (see also \cite{KGB} for an
alternative proof and generalization) $\P_{G(O)}(\Gr)\subset D_{G(O)}(\Gr)$
is a monoidal
subcategory. Moreover, $\P_{G(O)}(\Gr)$ is equipped with a commutativity
constraint
and a fiber functor, and we have
(for $k$ algebraically closed) 
an equivalence of Tannakian categories
$\P_\K(\Gr)\cong Rep(\LG)$.
This Theorem is known as the {\it geometric Satake isomorphism};
see \cite{doistor}, \cite{G}, \cite{MV} and \cite{BD}.
As the name suggests,
this result is a geometric, or categorical, counterpart of the classical
Satake isomorphism $K(Rep(\LG))\cong H_{sph}$,
where $H_{sph}$ is the spherical Hecke
algebra, and $K$ stands for the Grothendieck group.
Here the
word ``geometric'' means that, following the
 Grothendieck ``sheaf-function'' correspondence principle, one replaces
the space of functions on the set of $\Fq$-points of a scheme by the
category of $l$-adic complexes (or perverse sheaves) on this scheme
(or on its base change to an algebraically closed field).

In this and subsequent paper we extend the geometric Satake isomorphism
 to a description
of various categories of $l$-adic sheaves on $\Fl$ in terms of $\LG$.
These results have found several applications to representation theory:
to cells in affine Weyl groups and  bases in the
 Grothendieck groups of equivariant coherent sheaves \cite{nilkone};
to cohomology of tilting modules over quantum groups
at a root of unity \cite{cohotil}; and also
to Lusztig's conjectures on
 nonrestricted representations of  modular representations of $\Lg$
(in preparation, see announcement in \cite{ICM}).
We also think that they are closely related to some aspects of the recent work
\cite{GW} which discusses tamely ramified geometric Langlands duality from the point of
view of Yang-Mills theory.

The possibility to realize the affine Hecke algebra
$H$
and the ``anti-spherical'' module over it as Grothendieck
groups of (equivariant) coherent sheaves on varieties related to $\LG$
plays a crucial role in the proof of classification of irreducible
representations of $H$, which a particular case of the local Langlands conjecture,
  \cite{KL}, see also \cite{CG}. Thus one may hope
that the ``categorification'' of these realizations proposed here can 
contribute to the geometric Langlands program. 
Let us point out that existence of (some variant of) such a categorification
was proposed as a conjecture by V.~Ginzburg 
(see Introduction to \cite{CG}).

\subsubsection{}
Let us now describe some known statements about spaces of functions
on $G(F)$, whose geometric counterparts will be provided in the paper.

 Set $k=\Fpn$, and let $H=\Ce[I\backslash
G(F)/I]$ be the Iwahori-Matsumoto Hecke algebra. Let
 $T_w$ be the standard  basis
of $H$; here $w$ runs over the extended affine Weyl group $W$.
 Let $\Lambda\subset W$ be the coweight lattice of $G$, and
$\Lambda^+\subset \Lambda$ be the semigroup of dominant coweights. Let
$A\subset H$ be the commutative subalgebra generated by the elements
$T_\lambda$, $\lambda\in \Lambda^+$ and their inverses (see e.g.
\cite{doistor}, beginning of \S 7).
 Thus
$A$ has a basis $\theta_\lambda$, $\lambda\in \Lambda$, where
$\theta_\lambda$ are defined by the conditions $\theta_\lambda=
q^{-\ell(\lambda)/2}T_\lambda$
for $\lambda\in \Lambda^+$, $\theta_{\lambda+\mu}=\theta_\lambda\cdot
\theta_\mu$.

Recall that
 the {\it anti-spherical} (right\footnote{$H$ has a canonical anti-involution coming
 from the map $g\mapsto g^{-1}$, $g\in G(F)$; thus the categories of left
 and right modules are canonically identified. We define $M_{asp}$ as a right $H$-module
 to make some notations more natural: $M_{asp}$ is realized in the space of functions
 on $G(F)/I$ where $H$ acts naturally on the right.}) module
   $M_{asp}$ over $H$ is defined
as the induction
from the sign representation of the finite Hecke algebra $H_f\subset H$.
One can also describe $M_{asp}$ as follows. Let
 $C_w$ be the Kazhdan Lusztig basis of $H$;
let $W_f\subset W$ be the finite Weyl group, and $^fW\subset W$, 
$\fWf\subset W$ be the set
of minimal length representatives of, respectively, left and two-sided
 cosets of $W_f$ in $W$.
Then
\begin{equation}\label{Masp}
M_{asp}\cong H/\< C_w\ ,\ w\not \in ^fW\> .
\end{equation}
Notice that $M_{asp}$ is free of rank 1 over the subalgebra $A$.

Another important realization of $M_{asp}$ is in terms of the Whittaker model.
 Let $N\subset G$ be a maximal unipotent, and $\Psi: N(F)\to \Ce$
be a generic character. Then
\begin{equation}\label{MWh}
M_{asp}\cong (ind_{N(F)}^{G(F)}(\Psi))^I;
\end{equation}
here the right hand side is identified with the space of {\it Whittaker
functions on $G(F)/I$}.

The group $N(F)$ is not compact; because of this there is no straightforward
 definition of the category of Whittaker sheaves on $\Fl$
(the geometric counterpart of the right hand side
of \eqref{MWh}). Following \cite{GFV}
one can provide such a definition using Drinfeld's compactification
of the moduli space of $B$-bundles on a curve. However,
the following technically simpler (though probably less suited for
 generalizations)
approach  suffices for our purposes.

 Let $I_u\subset G(F)$
be the pro-$p$ radical of an Iwahori subgroup, %(possibly different
%from $I$), 
and $\psi: I_u\to \Ce$ be a generic character (the
definition is recalled below). Then one can use
 Lemma \ref{fW}
below to show that
\begin{equation}\label{MIWh}
M_{asp}\cong \left(ind_{I_u}^{G(F)}(\psi)\right)^I;
\end{equation}
and moreover,
the arising isomorphism between the right hand sides of
 \eqref{MWh} and \eqref{MIWh} is compatible with the standard bases
consisting of functions supported on one two-sided coset.

We call the right hand side of \eqref{MIWh} the Iwahori-Whittaker module.
It is easy to define the  category of Iwahori-Whittaker sheaves on $\Fl$.
 It  can probably
be shown to be equivalent to the category of Whittaker sheaves on
$\Fl$ (where the latter is defined following \cite{GFV}); this is not
pursued in this paper (see, however, Theorem \ref{Whitint} below).

\medskip

The methods of this paper can be used also to describe in a similar fashion
 geometric counterparts of the algebra $H$; this will be addressed in a future publication
 (see announcement in \cite{ICM}).

\medskip

Below we will define a triangulated monoidal category $D(\A)$
which is a geometric counterpart of the commutative algebra $A$; and
abelian categories $\fP$, $\PIW$
which are geometric counterparts of the right
hand sides of \eqref{Masp}, \eqref{MIWh} respectively;
%(see concerning the counterpart of \eqref{MWh});
we will then describe
$D(\A)$, and the derived categories $\Db(\fP)$, $\Db(\PIW)$,
  in terms
of the Langlands dual group.

\subsubsection{}
We now recall the realizations of $H$, $A$, $M_{asp}$ in terms of $\LG$,
whose categorical counterparts will be given below.

Let $\Lg$ be the Lie algebra of $\LG$. Let $\B=\LG/\LB$
be the flag variety; $\N$ be the nilpotent cone of $\Lg$,
and $p_{Spr}:\Nt=T^*(\LG/\LB)\to \N$ be the Springer map.
Let $St=\Ntil\times _\N\Ntil$ be the Steinberg variety of triples.

For a scheme  $X$ equipped with an action of an algebraic group $H$
we will let $Coh^H(X)$ be the category of $H$-equivariant  coherent sheaves;
we will write $D^H(X)$, or $D(X)$ if the group is unambiguous,
for the bounded derived category $ \Db(Coh^H(X))$. We will denote the Grothendieck 
group of either abelian or triangulated category $\C$ by $K(\C)$.

Let  $\bbH$ be the affine Hecke algebra of $G$.
 Thus $\bbH$  is an algebra over
$\Zet[v,v^{-1}]$, and $H=\bbH\otimes _{\Zet[v,v^{-1}]} \Ce$,
where the map $\Zet[v,v^{-1}]\to \Ce$ sends
$v$ to $q^{1/2}$.
One can define a subalgebra $\bbA\subset \bbH$, and a module $\bbM_{asp}$
such that $A=\bbA\otimes _{\Zet[v,v^{-1}]} \Ce$,
$M_{asp}=\bbM_{asp}\otimes _{\Zet[v,v^{-1}]} \Ce$.

Then we have an isomorphism (see e.g. \cite{CG} or \cite{Kthone})
$$\bbH\cong K(D^{\LG\times \Gm}(St));$$
where the algebra structure on the right hand side is
provided by  convolution.
Moreover, under this isomorphism
the subalgebra $\bbA$ is identified with the image of
$\delta_*: K(D^{\LG\times \Gm}(\Nt))\to  K(D^{\LG\times \Gm}(St))$,
where $\delta:\Ntil\imbed St$ is the diagonal embedding; notice
that $\delta_*$ is a homomorphism where the algebra structure
on $K(D^{\LG\times \Gm}(\Nt))$ is defined by $[\F]\cdot [\G]=[\F
\Lotimes _\O \G]$. The module  $\bbM_{asp}$ is identified
with  the $K(D^{\LG\times \Gm}(St))$ module  $K(D^{\LG\times \Gm}(\Nt))$.

\subsubsection{Informal summary}
Our method relies heavily on \cite{KGB} which provides a categorical counterpart of the description
of the center $Z(H)$ of the affine Hecke algebra $H$. According to a well known result of Bernstein
 \cite{doistor} we have $Z(H)\cong H_{sph}$, thus by Satake isomorphism we have $Z(H)\cong
K(Rep(\LG))$. In \cite{KGB} Gaitsgory uses geometric Satake isomorphism and nearby cycles functor
to define  a central functor $\Z$ from $Rep(\LG)$ to $D_I(\Fl)$ (the notion of a central functor is recalled
below). 

The present paper can be informally summarized as follows. We upgrade Gaitsgory's functor
$\Z$ to a functor from $D^{\LG}(\Nt)$, which is then shown to induce an equivalence with the 
Iwahori-Whittaker category, by linking various ingredients in the definition of $D^\LG(\Nt)$
to relevant structures on the perverse sheaves side. To make this more precise recall that
$\Nt=\{(\b, x)\ |\ \b \in \B, \ x\in rad(\b)\}$, where $\B$ is identified with the set of Borel subalgebras
in $\Lg$ and $rad$ stands for the nilpotent radical. We show that, in the appropriate formal sense,
the tensor functor from $Rep(\LG)$ to $D_I(\Fl)$ corresponds to the fact that in the dual side
we deal with the category of $\LG$ equivariant coherent sheaves on some algebraic variety.
The element $x\in \Lg$ in the description of $\Nt$ arises from the logarithm of monodromy acting
on the nearby cycles sheaf $\Z(V)$, $V\in Rep(\LG)$ by Tannakian formalism. 
Finally, the "flag" $\b\in \B$ corresponds to a filtration on $\Z(V)$, $V\in Rep(\LG)$ by
{\it Wakimoto sheaves}, see Theorem \ref{filtrWakr}.
Wakimoto sheaves categorify elements $\theta_\la\in A\subset H$, and 
 Theorem \ref{filtrWakr} is a categorification of the fact that $Z(H)\subset A$, thus an element in 
 $Z(H)$ is a linear combination of $\theta_\la$, $\la\in \La$.
 On the other hand,  Theorem \ref{filtrWakr} is
equivalent to the computation of cohomology of the
so-called {\it semi-infinite} orbits with coefficients in these sheaves,
and is closely related to Mirkovi\' c and Vilonen's computation of corresponding
cohomology for spherical sheaves;
see section \ref{filtrWak} for further comments.

We finish the introduction by pointing out another result on the structure of  central
sheaves of \cite{KGB}, Theorem \ref{tilt}
proved below. 
It says that the objects of the 
Iwahori-Whittaker  category cooked out of central sheaves are {\it tilting}.
This result is inspired by the  {\it ``Koszul duality''} yoga
of \cite{BGS}, see Remark \ref{Koszul}.

\medskip

Finally, let us make a standard remark  that all the  
 results and proofs of the paper
work in the alternative setup where the finite characteristic base field
 $k$ is replaced by $\Ce$, the field of coefficients
$\Ql$  is also replaced
by $\Ce$, and the category of $l$-adic constructible sheaves by
 the category of $D$-modules (in the part of the paper where neither
Artin-Schreier sheaf, nor weights are used one can work with 
 constructible sheaves in the classical topology).

\subsection{More notations}\label{nota}
From now on we fix $k=\Fpbar$.

The convolution diagram will be written as $\Fl\underset{\bI}{\times}
\Fl\to \Fl$. If $X\subset \Fl$, $Y\subset \Fl$ are subschemes, and
$Y$ is $\bI$ invariant, then we get a subscheme $X\underset{\bI}{\times}
Y\subset \Fl\underset{\bI}{\times}\Fl$. For $\F\in D(\Fl)$, $\G\in
D_I(\Fl)$ we get an object $\F\underset{\bI}{\boxtimes}\G$
(twisted product) of the category of $l$-adic complexes on
$\Fl\underset{\bI}{\times}\Fl$.

Let $\kappa:\Lambda\to \fW$ be the bijection such that $\kappa(\lambda)
\in  W_f\cdot \lambda$.  

All derived categories below will be bounded derived categories, notation $D$ will be used
instead of a more traditional $D^b$, unless stated otherwise.

We now describe the results of the paper.

\subsection{The monoidal functor}
If $X$ is smooth, then $D(X)$ is a tensor category under the
(derived) tensor product of coherent sheaves. The first result of
the paper (see section \ref{constrF}) is construction of a monoidal
functor
\begin{equation}\label{FtoDI}
%\begin{array}{ll}
D^\LG(\Ntil)\to D_I. 
%\\
%D^{\LG\times \Gm}(\Ntil)\to D_I^{mix}.
%\end{array}
\end{equation}
In fact we will do a little bit more. We will define (in
section \ref{A})  a full 
subcategory $\A\subset \P_I$ %, $\A^{mix}\subset \P_I^{mix}$
 which will turn out to be closed under 
convolution. Then the homotopy category $Hot(\A)$,
%  $Hot(\A^{mix})$ 
of finite complexes of
objects in $\A$
inherits a monoidal structure.
 Further, let $D(\A)$ %, $D(\A^{mix})$
 denote the quotient of the triangulated category %ies
$Hot(\A)$ %, $Hot(\A^{mix})$ respectively
 by the subcategory of acyclic complexes.
Then $D(\A)$ %,  $D(\A^{mix})$ are 
is also a monoidal category.
 We have  the obvious functor $D(\A)
\to \Db(\P_I)$.
%,  $D(\A) \to \Db(\fP_I^{mix})$.

We will construct a monoidal functor
\begin{equation}\label{F}
%\begin{array}{ll}
F:D^\LG(\Ntil)\to D(\A)
%\Fmix:D^{\LG\times \Gm}(\Ntil)\to D(\A^{mix}).
%\end{array}
\end{equation}
%here the monoidal structure on $D^\LG(\Ntil)$, 
%$D^{\LG\times \Gm}(\Nt)$ is provided by the (derived) tensor product.

%We then describe some natural module categories for the monoidal
%categories $D(\A)$, $D(\A^{mix})$
% such that for a certain object $X$ of the category
%the functor $\F\mapsto F(\F)(X)$ from $D(\Nt)$ to the category in
%question is an equivalence.

One can use the argument of \cite{Be}, 1.3 to define a natural
functor $\Db(\PI)\to D_I$; thus we can define \eqref{FtoDI}
as the composition $D^\LG(\Nt)\to D(\A)\to \Db(\P_I)\to D_I$.
 It is then easy
to see that it comes with a natural monoidal structure.

Below we will not use \eqref{FtoDI}, only  \eqref{F}.

\subsection{Compatibility with Frobenius}
Let $\bq: \Nt\to \Nt$ be the map sending a pair $(b,x)\in \Nt\subset
\B \times \Lg$ to $(b,qx)$. %%or q^{-1}???

Let also $Fr=Fr_q$ be the geometric Frobenius; 
recall that  for a (ind)scheme $X$ over $\Fq$, $Fr$ induces an autoequivalence
of the (derived) category of $l$-adic sheaves on  $X_{\Fpbar}$,
$X\mapsto Fr^*(X)$, see \cite{Weil2}.

\begin{Prop}\label{Frq}
There exists a natural isomorphism of functors
\begin{equation}\label{Frqeq}
Fr^*\circ F \cong F\circ \bq^*
\end{equation}
\end{Prop}

\begin{Rem} Let the multiplicative group $\Gm$ act on $\Nt$
by   $t:(b,x)\mapsto (b,t^{-2}x)$. Then for $\F,\G\in D^{\LG\times \Gm}
(\Nt)$ the vector space $Hom_{D^\LG(\Nt)}(\F,\G)$ carries an
action of $\Gm$, hence a grading. The Proposition provides isomorphisms
$Fr^*(F(\F)) \cong \F$, $Fr^*(F(\G))\cong \G$;
the map $Hom_{D^\LG(\Nt)}(\F,\G)\overset{F}{\To} Hom_{D(\Fl)}(F(\F),F(\G))$
carries the degree $n$ component into a subspace where  Frobenius 
acts with weight $n$.

It follows that the equivalences $\fPhi$, $\FIW$ defined below also satisfy
this property.
\end{Rem}
\begin{Rem}
In fact, a slight modification of our argument provides a monoidal functor
from $D^{\LG\times \Gm}(\Nt)$ to the derived category of mixed
$l$-adic sheaves on $\Fl_{\Fq}$. We expect that it induces an equivalence
between $D^{\LG\times \Gm}(\Nt)$ and $D(\fP_I^{mix})$, $\DIW^{mix}$,
where  $D(\fP_I^{mix})$, $\DIW^{mix}$ are {\it mixed versions} of the
categories  $D(\fP)$, $\DIW$ in the sense of \cite{BGS}.
(We warn the reader that the mixed category is {\it not} the category
of all mixed sheaves with an appropriate equivariance condition;
a necessary condition for a perverse sheaf to lie in the mixed category
is that its  associated graded
with respect to the weight filtration is semisimple,\footnote{We thank
V.~Ginzburg for pointing out this difficulty to us.} cf \cite{BGS}, 4.4). 
\end{Rem}

\subsection{Anti-spherical quotient category}
Recall that $\bI$ orbits on $\Fl$ (the so-called Schubert cells)
are parameterized by $W$; for $w\in W$ let $j_w:\Fl_w\imbed \Fl$
be the embedding of the corresponding Schubert cell. We let
 $L_w=j_{w!*}(\cons[\ell(w)])$,
$w\in W$ be the irreducible objects of $\P_I$, and
 $j_{w!}=j_{w!}(\cons[\ell(w)])$,  $j_{w*}=j_{w*}(\cons[\ell(w)])$
be the standard and costandard objects.
%; for finite $k$  (the ``mixed'' case) we set
% $L_w=j_{w!*}(\cons[\ell(w)])(\frac{\ell(w)}{2})$,
% $j_{w!}=j_{w!}(\cons[\ell(w)])(\frac{\ell(w)}{2})$,
%  $j_{w*}=j_{w*}(\cons[\ell(w)])(\frac{\ell(w)}{2})$.
For an abelian category $\A$, and a set $S$ of irreducible
objects of $\A$ let $\<S\rangle$
denote the full abelian subcategory of objects obtained from elements
of $S$ by extensions.
Define the  Serre quotient category %ies
of $\PI$ %, $\PI^{mix}$
 by
$$\fP_I=\P_I/\< L_w\ |\ w\not \in \fW\>.$$
%$$\fP_I^{mix}=\P_I^{mix}/\< L_w(i)\ |\ w\not \in \fW, \ i\in \Zet\>.$$

%%%%%%%%%%%%%%mixed from here 

Let $pr_f :\P_I\to \fP_I$ %, $pr_f^{mix}:\P_I^{gr}\to \fP_I^{mix}$
 be the projection functor.

\begin{Thm}\label{Phi}
%a)
 The functor $\fPhi:=pr_f\circ F$ is an equivalence
\begin{equation}\label{feq} %\begin{array}{ll}
\fPhi:D^\LG(\Ntil)\iso D(\fP_I).
\end{equation}

%b) For $\F,\G\in D^{\LG\times \Gm}(\Nt)$ the isomorphism
%\begin{equation}\label{isom}
%Hom_{D^\LG}(\F,\G)\iso Hom _\Fqbar(\fPhi^{mix}(\F),
%\fPhi^{mix}(\G))
%\end{equation}
%is compatible with the grading. Here
% $\fPhi^{mix}$ :=pr_f^{gr}\circ F^{gr}$; the grading
%on the left hand side of \eqref{isom} is provided by the $\Gm$
%action, and the grading on the right-hand side of \eqref{isom}
%is by Frobenius weights. 
\end{Thm}

\subsection{Iwahori-Whittaker category}
Let $B=T\cdot N$, $B_-=T\cdot N_-$ be  opposite Borel subgroups, and
assume that $\bI, \bI_-\subset \GO$,\  $\bI\supset \BO$, \
 $\bI_-\supset \BOm$ for  Iwahori group schemes $\bI$, $\bI^-$.
Let $\bI_u^-\subset \bI^-$ be the pro-unipotent radical.

Let also $\NFm\subset \GF$ be the group ind-scheme, $\NFm(k)=N^-(F)$.
For a simple root $\alpha$ let $u_\alpha:\NFm\to \Ga$ be the
corresponding homomorphism. We define $\Psi:\NFm\to \Ga$ by
$\Psi(n)=Res(\frac{dt}{t} \sum u_\alpha(n))$. We also define
$\psi_I:\bI_u^-\to \Ga$ by $\psi_I(g^-g^{\geq 0} )= \Psi(g^-)$ for
$g^-\in \bI_u^-\cap \NFm$, $g^{\geq 0}\in \bI_u^-\cap \BF$. Let
$\DIW$, $\PIW$ be, respectively, the $\bI^u_-,\psi$-equivariant derived category
of $l$-adic sheaves on $\Fl$, and the subcategory of perverse sheaves therein.
Since the group scheme $\bI^u_-$ is pro-unipotent, it follows that
the forgetful functor $\DIW$ to $D(\Fl)$ is a full embedding, thus $\DIW\subset D(\Fl)$,
$\PIW\subset \P(\Fl)$ are full subcategories. 
 
Let us mention the following
\begin{Lem}\label{fromBGS}
We have a natural equivalence
$$\DIW\cong \Db(\PIW).$$
%$$\DIW^{gr}\cong \Db(\PIW^{gr}).$$
\end{Lem}
\proof is parallel to that of Lemma 4.4.6, and Corollary 3.3.2
in \cite{BGS}. \epf

We have an injection $W\imbed \Fl$, $w\mapsto w\bI$; for
each one of the (ind) group schemes
 $\NFm$,  $\bI_-^u$, $\bI$ the image of this
map is a set of representatives for the orbits of the group on
$\Fl$. For $w\in W$ let $\Fl^w$ (respectively, $\Sigma_w$)
 be the corresponding orbit of
$\bI_-^u$ (respectively, $\NFm$),
 and $i_w:\Fl^w\imbed \Fl$,  $\iota_w:\Sigma_w\imbed \Fl$ be the embeddings.

The proof of the next Lemma is left to the reader.

\begin{Lem}\label{fW} a) For $w\in W$ the following are equivalent

i) $w\in ^fW$.

ii) $Stab_{\bI^-_u}(w\bI) \subset Ker(\psi_I)$.

iii)  $Stab_{\NFm}(w\bI)\subset Ker (\Psi)$.

b) For $w$ satisfying the equivalent conditions in (a) 
the $\bI^-$-orbit $\Fl^w$ is contained in one $\NFm$-orbit. \epf
\end{Lem}

 If $w\in \fW$ then there exist unique
maps $\psi_w:\Fl_w\to \Ga$, $\Psi_w:\Sigma_w\to \Ga$ defined by
$\psi_w(g\cdot w\bI)=\psi(g)$, $\Psi_w(n\cdot w\bI)=\Psi(n)$.
Define  $\Delta_w, \nabla_w\in \PIW$ by
$\Delta_w=i_{w!}\psi_w^*(\AS)[\ell(w)]$,
$\nabla_w=i_{w*}\psi_w^*(\AS)[\ell(w)]$ where $\AS$ in the
Artin-Schreier sheaf.

Define the functor $Av_\Psi: D_I\to \DIW$,
%$Av_\Psi^{gr}: D_I^{gr}\to \DIW^{gr}$
 by $\F\mapsto \Delta_0*\F$.

\begin{Thm}\label{IWisfP}
%a) We have $$Av_\Psi:\P_I \mapsto \PIW;$$
%and $Av_\Psi|_{\P_I}$ factors through $\fP_I$.
%b)
The functor  $Av_\Psi|_{\P_I}$ %,  $Av_\Psi^{gr}|_{\P_I^{gr}}$
 induces an equivalence $\fP_I\iso \PIW$.
% $\fP_I^{gr}\iso \PIW^{gr}$.
\end{Thm}

Define  the functor $\FIW:D^\LG(\Nt)\to \Db(\PIW)=\DIW$
%, $\FIW^{gr}:D^{\LG\times\Gm}(\Nt)\to \Db(\PIW^{gr})=\DIW ^{gr}$ 
by $\F\mapsto Av_\Psi \circ F(\F)$.

In view of Theorem \ref{IWisfP} and Lemma \ref{fromBGS}, Theorem \ref{Phi} is
equivalent to the following

\begin{Thm}\label{descrIWiprec}
 The functor $\FIW$ provides an equivalence
$D^\LG(\Nt)\cong \DIW$.

\end{Thm}

\end{section}

\begin{section}{Comparison of anti-spherical and Whittaker categories}
In this section we will prove a result in the direction
of Theorem \ref{IWisfP}. The proof of the Theorem will be finished in
section \ref{IWisfPkonec} after the proof of Theorem \ref{descrIWiprec}.
%In this section we skip the evident ``mixed'' variant of the statements.

\begin{Prop}\label{IWisfPprop}
a) We have $$Av_\Psi:\P_I \mapsto \PIW;$$
thus $Av_\Psi|_{\P_I}$ induces an exact functor $\P_I\to \PIW$.
This functor factors through $\fP_I$.

b) The functor $\fP_I\to \PIW$ induced by $Av_\Psi$ is a full embedding.
\end{Prop}

Set $\delta_e=j_{e*}=j_{e!}$ where $e\in W$ is the identity element.

Let $W'\subset W$ be the subgroup generated by simple
reflections (non-extended affine Weyl group). Thus $\cupl_{w\in W'}\Fl_w$ 
is a connected component of $\Fl$.

\begin{Lem}\label{deltae}
a) For $w\in W'$ we have nonzero
morphisms % in $D_I$
$$\delta_e\to j_{w!};$$
$$j_{w*}\to \delta_e,$$
whose (co)kernel does not contain $\delta_e$ in its Jordan-Hoelder
series.

b) If $w=w_1w_2\in W$, $w_2\in W'$  and $\ell(w)=\ell(w_1)+\ell(w_2)$ then
$$\dim Hom(j_{w_1!}, j_{w!})=1=\dim Hom(j_{w*}, j_{w_1*});$$
and a nonzero map $j_{w_1!}\to j_{w!}$ (respectively,
$j_{w*}\to j_{w_1*}$) is injective (respectively, surjective).
\end{Lem}
\proofpt We prove the statements concerning $j_{w!}$, the ones concerning
$j_{w*}$ are obtained by duality.

For a simple reflection $s_\alpha\in W$
we have an exact sequence of perverse sheaves on the projective
line $\overline{\Fl_{s_\alpha}}$
\begin{equation}\label{ontheline}
0\to  \delta_e\to j_{s!} \to L_s\to 0.
\end{equation}
If $u\in W$ is such that $\ell(u\cdot s_\alpha)>\ell(u)$ consider the
convolution of $j_{u!}$ with \eqref{ontheline}; it is an exact triangle
\begin{equation}\label{ontheotherline}
 j_{u!}\to j_{u\cdot s!} \to j_{u!}* L_{s_\alpha}.
\end{equation}
Notice that $j_{u!}* L_s=\pi_\alpha^*\pi_{\alpha*}(j_{u!})[1]$,
 where $\pi_\alpha:\Fl\to \Fl(\alpha)$
is the projection to the partial affine flag variety $\Fl(\alpha)
=\GF/\bI_\alpha$ for the minimal parahoric $\bI_\alpha$ corresponding to
$\alpha$. Since $\pi_\alpha\circ j_u$ is a locally closed affine embedding
(because $\ell(u\cdot s_\alpha)>\ell(u)$), we see that
$\pi_{\alpha*}(j_{w!})$, and hence $\pi_\alpha^*\pi_{\alpha*}(j_{u!})[1]$
are perverse sheaves. Thus the exact triangle \eqref{ontheotherline}
is in fact an exact sequence of perverse sheaves.
 Also all irreducible
subquotients of $j_{u!}* L_{s_\alpha}$ are of the form $\pi_\alpha^*(L[1])$
for a perverse sheaf $L$ on $\Fl(\alpha)$; thus none of them is isomorphic
to $\delta_e$. This implies (a) by induction in $\ell(w)$.

Since $j_{w_2!}$ is invertible under convolution (see Lemma
\ref{Lemmult}(b) below) we have
$$Hom(j_{w_1!}, j_{w!})=Hom(j_{w_1!}, j_{w_1!}*j_{w_2!})=
Hom(\delta_e,j_{w_2!}),$$
thus the first statement in (b) follows from (a). Finally, the exact
sequence \eqref{ontheotherline} implies by induction in $\ell(w_2)$
existence of an injective map $j_{w_1!}\to j_{w!}$. \epf

\begin{Lem}\label{nunu}
a) We have %\begin{equation}\label{Avnul}
$$Av_\Psi(L_w)=0\ \ {\iff} \ \ w\not \in ^fW.
$$
%\end{equation}

b) We have $\Delta_e\cong \nabla_e$.

c) For   $w=w_f\cdot w'$, $w_f\in W_f$, $w'\in \fW$ we have
%\begin{equation}\label{Avst}
$$\begin{array}{ll}
&Av_\Psi(j_{w!})\cong \Delta_{w'},\\
 &Av_\Psi(j_{w*})\cong \nabla_{w'}.
\end{array}
$$
%\end{equation}

\end{Lem}
\proofpt If $w\in \fW$, then the convolution map $\Fl\underset{\bI}{\times}
\Fl\to \Fl$ restricted to the generic point of the support of $\Delta_0
\underset{\bI}{\boxtimes}L_w$ is an isomorphism. Hence $\Delta_0*L_w
\ne 0$ for $w\in \fW$. On the other hand,
for $w\not \in \fW$ there exists a simple root $\alpha\ne
\alpha_0$ such that $L_w$ is equivariant with respect to the  corresponding
minimal parahoric subgroup $\bI_\alpha$ (here $\alpha_0$ denotes the affine
simple root). Then the functor $\F\mapsto \F*L_w$ factors through
the functor $\pi_{\alpha*}$ (recall that
$\pi_\alpha:\Fl\to \Fl(\alpha)$ is the projection to the corresponding partial
affine flag variety). However, $\pi_{\alpha*}(\Delta_0)=0$ because
the character $\psi_I$ is nontrivial on $Stab_{\bI_u^-}(x)$ for
any $x$ in the image of the support of $\Delta_0$ under $\pi_\alpha$.
This proves (a).

(b) is clear because  $\psi_I$ is nontrivial on $Stab_{\bI_u^-}(x)$ for
any $x\in \overline{\Fl^e}-\Fl^e$.

In view of (b) it suffices to prove
 the first equality in (c); the second one then follows by duality.
The equality is clear when $w\in \fW$, because in this
case the convolution map restricted to the support of $\Delta_0
\underset{\bI}{\boxtimes}j_{w!}$
 is an isomorphism over $\Fl^w$, while the $*$
restriction of $\Delta_0
\underset{\bI}{\boxtimes}j_{w!}$ to the preimage of the complement of
$\Fl^w$ is zero. Let now $w$ be arbitrary; we have $w=w_f
\cdot w'$ for some 
$w_f\in W_f$, $w'\in \fW$, where $\ell(w)=\ell(w_f)+ \ell(w')$.
Then Lemma \ref{deltae} and part (a) of this Lemma imply
that $\Delta_e*j_{w_f!}\cong \Delta_e*\delta_e=\Delta_e$. Thus
we have $$\Delta_e*j_{w!}\cong \Delta_e*j_{w_f!}*j_{w'!}\cong
\Delta_{w'}.\ \ \ \epf$$

For an algebraic group
$H$ and a subgroup $H'\subset H$ (or more generally, for group
schemes of possibly infinite type, such that the quotient $H/H'$
is of finite type) let $\Gamma_{H'}^H$ be the $*$ induction
functor from $H'$-equivariant to $H$-equivariant sheaves; recall
that it is defined by $\Gamma_{H'}^H(\F)=
a_*(\cons\boxtimes_{H'}\F)$, where $a:H\times_{H'}X\to X$ is the
action map (cf \cite{BL}).

Define the functor $Av_{I}:\DIW\to D_I$ by   $Av_{I}=\Gamma_{\bI\cap
\bI_u^-}^\bI$. % (notations of section \ref{gammatam} below).

\begin{Lem}\label{Hpo}
 We have $H^{p,0}(Av_{I}(\Delta_0)[\ell(w_0)])
\cong j_{w_0!}$, where $w_0\in W_f$
is the longest element, and superscript $^p$ refers to the $t$-structure of perverse sheaves.
% $pr_f(j_{w_0*})\cong pr_f(j_e)$ is irreducible.
\end{Lem}

\proofpt It suffices to construct an exact triangle
$$  j_{w_0!}\to Av_I(\Delta_0)[\ell(w_0)]\to C$$
such that $C\in D^{p,>0}$. The definition of $Av_I$ implies
that $$Hom_{D_I}(X,Av_I(\Delta_0)[\ell(w_0)])=Hom_{D}
(Forg(X),(\Delta_0)[\ell(w_0)]),$$
where $D$ is the derived category of $l$-adic sheaves on $\Fl$,
and $Forg:D_I\to D$ is the forgetful functor (notice that
$D_{\bI\cap \bI_u^-}$ is a full subcategory in $D$ because $\bI\cap \bI_u^-$
is unipotent). The proof of Lemma \ref{nunu}(a) shows that
\begin{equation}\label{LwD}
Hom^\bu(L_w,\Delta_0)=0
\end{equation}
for $w\in W_f$, $w\ne e$ (and more generally for $w\not \in \fWf$).
Also it is clear that $$Hom_{D_I}(\delta_e,
Av_I(\Delta_0)[\ell(w_0)])\cong \Ql.$$
Now Lemma \ref{deltae}(a) implies that
$$Hom_{D_I}(j_{w!}, Av_I(\Delta_0)[\ell(w_0)])\cong \Ql$$
for $w\in W_f$.
Moreover,
 the composition of nonzero arrows
$$j_{w!}\to j_{w'!}\to Av_I(\Delta_0)[\ell(w_0)]$$
is nonzero, because the composition $\delta_e\to j_{w!}\to j_{w'!}$
is nonzero by \ref{deltae}(b). 
 Hence, if $C=cone(j_{w_0!}\to Av_I(\Delta_0)[\ell(w_0)])$,
then using again Lemma \ref{deltae}(b) we see that for all $w\in W_f$
$$Hom(j_{w!},C[i])=0$$
for $i\leq 0$, which implies that $C\in D^{p,>0}$ according to the
definition of perverse $t$-structure. \epf

\begin{Rem}\footnote{This Remark is included here following the referee's suggestion.}
A slightly different proof of the Lemma can be given as follows.
It is not hard to describe  $Av_{I_u}(\Delta_0)$, 
where $\bI_u\subset \bI$ is the pro-unipotent
radical (details will appear in \cite{BM}).
  This is obviously an object in the category of $\bI_u$ equivariant sheaves
supported on $G/B\subset \Fl$, which is identified with
category $O$ for $G$. Then $Av_{I_u}(\Delta_0)\cong \Xi [\ell(w_0)]$,
 where $\Xi$ is the maximal
projective in category $O$ 
(a projective cover of the irreducible Verma module).
Lemma \ref{Hpo} can then be deduced from the fact that the subobject 
$j_{w_0!}\subset \Xi$ 
is the maximal $B$-equivariant subobject in the $B$-monodromic object $\Xi$.
\end{Rem}

\subsubsection{Right inverse to $Av_\Psi$} To prove Proposition
 \ref{IWisfPprop} we will explicitly construct a right inverse functor to $Av_\Psi$.
Namely, define $F':\PIW\to \fP_I$ by
 $$F'(\F)= pr_f( H^{p,\ell(w_0)}( Av_I(\F))).$$
To motivate this definition we remark that one can easily show that
$F'$ is right adjoint to $Av_\Psi$ (we neither check, nor use this fact below).
\begin{Lem}\label{rightin}
There exists a canonical  isomorphism
%\begin{equation}\label{FprF}
$F'\circ Av_\Psi\cong id.$
%\end{equation}
\end{Lem}
\proofpt
For $w\in W_f$ the functor from $D_I$ to $^fD_I$ sending $\F$
to $pr_f(L_w*\F)$ is zero if $w\ne e$ and is isomorphic to $pr_f$
otherwise. Hence the convolution functor descends to a functor
$\fP_I^0 \times \fP_I\to \fP_I$ exact in each variable; here
$\fP_I^0$ denotes the Serre quotient category of $\bI$-equivariant
perverse sheaves on $\GO/\bI\subset \Fl$ by the subcategory generated
by $L_w$, $w\ne e$. In particular, for $\F\in \P_I$ we have
$$ F'\circ Av_\Psi(pr_f(\F))\cong
 pr_f\circ H^{p,0}(Av_I(\Delta_0[\ell(w_0)])*\F)
\cong  pr_f\circ H^{p,0}(j_{w_0!}*\F)\cong pr_f(\F),$$
where the last isomorphism follows from Lemma \ref{deltae}(a), and
the previous one from Lemma \ref{Hpo}. \epf

\subsubsection{Proof of Proposition \ref{IWisfPprop}(conclusion)}\label{prooW}
For a triangulated category $D$ and a set of objects $\cS\subset
Ob(D)$ we let $\< \cS\rangle$ be the set of all objects obtained
from elements of $\cS$ by extensions; i.e. $\< \cS\rangle$ is the
smallest subset of $D$ containing $\cS\cup\{0\}$
 and such that
for all $A,B\in \<\cS\rangle$ and an exact triangle $A\to C\to
B\to A[1]$ we have $C\in \cS$.

The definition of perverse $t$-structure implies that
$$Ob(\P_I)
=\< j_{w!}[i]\ |\ i\geq 0\rangle \cap \< j_{w*}[i]\ |\ i\leq 0\rangle;$$
$$Ob(\PIW)
=\< \Delta_{w}[i]\ |\ i\geq 0\rangle \cap \< \nabla_{w}[i]\ |\ i\leq
0\rangle.$$
Thus the first statement in Proposition \ref{IWisfPprop}(a) follows
from Lemma \ref{nunu}(c). The second one is immediate from part
(a) of that Lemma.
Part (a) of the Proposition is proved.

We also see that $F(L_w)$ is irreducible for $w\in \fW$, because it is
the image of a nonzero map $\Delta_w\to \nabla_w$; it is clear that
$F(L_w)\not \cong F(L_{w'})$ for $w\ne w'$, $w,w'\in \fW$,
because they have different supports.

Thus the Proposition follows from Lemma \ref{rightin} and the following
Lemma. \epf

\begin{Lem}
Let $F:\A\to \B$ be an additive functor between
abelian categories. Assume that

i) $F$ is exact.

ii) Every object of $\A$ has finite length, and $F$ induces an isomorphism
$Hom(L_1,L_2)\iso Hom(F(L_1), F(L_2))$ for any irreducible objects
$L_1,L_2$ of $\A$.

iii) There exists an additive functor $F':\B\to \A$ such that
$F'\circ F\cong id$.

Then $F$ is a full embedding.
\end{Lem}

\proofpt Conditions (i) and (iii) imply that $F$ is injective on $Ext^1$.
Indeed, let  $0\to X\to Y\to Z\to 0$ be  a short exact sequence in $\A$.
If $F(Y)\cong F(X)\oplus F(Z)$ is the splitting of its image under $F$,
then applying $F'$ to it we see that the original sequence is split.

Now induction in the lengths of $X,Y$
shows that  $F$ induces an isomorphism $Hom(X,Y)\to
Hom(F(X),F(Y))$ for any two objects $X,Y\in \A$. \epf

\end{section}
\begin{section}{Construction of the monoidal functor
$F$}\label{constrF}
\subsection{Plan of the construction} We will use a version of Serre's
description of $Coh({\mathbb P}^n)$ as a quotient of the category of graded modules
over the symmetric algebra.
We need some notations.

Let $\Nthat$ be the preimage of
$\Nt\subset \B\times \Lg$ in $\LGU\times \Lg$. Thus $\Nthat$ is a
locally closed subscheme in the affine variety $\LGUbar\times\Lg$;
here $\LU\subset \LG$ is a maximal unipotent, and
$\LGUbar$ is the affine closure of the basic affine space
$\LG/\LU$. We now define a closed (obviously affine) 
subscheme  $\Nthat_{af}\subset \LGUbar\times\Lg
$ containing $\Nthat$ as an open subscheme (though {\em different} from the closure of
$\Nthat$ in  $\LGUbar\times\Lg$).
 On $\LGU\times \Lg $ we have a canonical vector field
$v_{taut}$ whose value at a point $(p,x)$ equals $(a(x), 0)$ where
$a$ stands for the action of the Lie algebra $\Lg$ on $\LGU$. The
vector field $v_{taut}$ induces a derivation of $\O_{\LG/\LU\times
\Lg}$; we let $\Nthat_{af}$ be the zero-set of this derivation
(i.e. the defining ideal of $\Nthat_{af}$ is generated by the
image of the derivation). It is clear that $\Nthat_{af}\cap
(\LGU\times \Lg)=\Nthat$. We set $\Ohat_{\Nt}=\O_{\Nthat_{af}}$,
and call this ring the multi-homogeneous coordinate ring of $\Nt$.
For a scheme $S$ equipped with an action of an algebraic group $H$
we write $Coh_{fr}^H(S)$ (or $\O-mod^H_{fr}$ is $S=Spec(\O)$) for
the full subcategory in $Coh^H(S)$ consisting of objects of the
form $V\otimes \O_S$, $V\in Rep(H)$.

We now describe the plan of the construction.
The first piece of data is a monoidal functor $\Fbar: Rep(\LG\times
\LT)\to \P_I$ (i.e. a monoidal functor to $D_I$
landing in $\P_I$); the main ingredient is provided  by \cite{KGB}.

We then explain that a certain natural endomorphism of this action
(also defined in \cite{KGB})
yields an extension of $\Fbar$ to a monoidal functor $\Ftil:
Coh_{fr}^{\LG\times \LT}(\LGUbar \times \Lg)\to \P_I$. For this
 we describe in section \ref{ar}
certain distinguished arrows in $Coh_{fr}^{\LG\times
\LT}(\LGUbar)$; and also  arrows between objects $\Fbar(V)$,
$V\in Rep(\LG\times \LT)$. We will then require 
the functor $\Ftil$ to intertwine the two sets of arrows.
Some formalism described in section
\ref{form} shows that this requirement defines $\Ftil$ uniquely
(part of the argument  is a variation of the standard description
of elements in a Lie algebra as tensor endomorphisms of a fiber
functor).

Then $\Ftil$
is constructed; it yields a functor (again denoted by
$\Ftil$) from $Hot(Coh^{\LG\times \LT}_{fr}(\Nthat_{af}))$ to
$Hot(\PI)$, where $Hot$ stands for the homotopy category.
Let $Acycl\subset Hot(Coh^{\LG\times \LT}_{fr}(\Nthat_{af}))$ be
the full subcategory of such complexes $\F^\bu$ that
$\F^\bu|_{\Nthat}$ is acyclic. In section \ref{filtrWak} we prove
certain facts about the central sheaves of \cite{KGB}, and deduce
from it that $\Ftil$ sends ${Acycl}$ to acyclic complexes.
Hence $\Ftil$
factors to a  functor $D^\LG(\Nt)\overset{F}{\To}
\Db(\PI)\to D_I$.

\subsection{Central and Wakimoto sheaves: definition of the functor
$\Fbar$}\label{mlm}
%\subsubsection{}
 Recall the functor $\Z:\P_\K(\Gr)\to \P_I\subset D_I$
constructed in \cite{KGB}.
 
%The results of {\it loc. cit.}
%also give a central functor
%$\F_\K^{mix}(\Gr)\to D_I^{mix}$, which we also denote by $\Z$.

We identify  $\P_\K(\Gr)$ with $Rep(\LG)$ by means of the geometric Satake
equivalence $\Sa:Rep(\LG)\to \P_\K(\Gr)$. 
We set 
$V_\lambda=\Sa^{-1}(IC_\lambda)$ where $IC_\lambda=\jbar_{\lambda !*}\left(
\cons[\ell(\lambda)]        %(\frac{\ell(\lambda)}{2})
\right)$,
 and
$\jbar_\lambda:\Gr_\lambda\imbed \Gr$ is the embedding
of  the image $\Gr_\lambda$ of $\Fl_\lambda$ under the projection
$\pi:\Fl\to \Gr$; thus $V_\lambda$ is a representation with highest weight
$\lambda$. 
Notice that the convolution map $supp(IC_\lambda\underset{\GO}{\boxtimes}
IC_\mu) \to supp(IC_\lambda*IC_\mu)=\overline{\Gr_{\lambda+\mu}}$
is an isomorphism over $\Gr_{\lambda+\mu}$; hence we have
$$\jbar_{\lambda+\mu}^*(IC_\lambda\underset{\GO}{\boxtimes}
IC_\mu)\cong \cons[\ell(\lambda+\mu)] %(\frac{\ell(\lambda)}{2})
$$
canonically. Thus we get a canonical element $m_{\lambda,\mu}$ in
 the one dimensional vector space $Hom(IC_\lambda*IC_\mu,
IC_{\lambda+\mu})$.

 We also set $Z_\lambda=\Z(V_\lambda)$.

%\subsubsection{}
 The functor $\Z$ is monoidal, and
moreover {\em central}; the latter means that for every $V\in Rep(\LG)$ and
 $\F\in D_I$ there is a fixed ``centrality''
isomorphism $\sigma_{V,\F}:\Z(V)*\F\cong 
\F*\Z(V)$ satisfying
some natural compatibilities (spelled out e.g. in \cite{B},
\S 2.1, and checked in \cite{KGB} and Gaitsgory's Appendix to \cite{B}).
 Notice that a central functor from a tensor category $\A$ to a monoidal
category $\C$ is the same as a tensor (compatible with braiding)
functor from $\A$ to the center of $\C$ (see e.g. \cite{Kassel}, XIII.4).

Recall that
$j_{w!}=j_{w!}(\cons[\ell(w)]) %(\lT) 
)$, $j_{w*}=j_{w*}(\cons[\ell(w)] %(\lT)
)$;
and  $\delta_e=j_{e!}=j_{e*}$ is the unit object of $D_I$ (here
$e$ is the unit element of $W$). The following statement is
well-known.
\begin{Lem}\label{Lemmult} a) If $w_1, w_2 \in W$ are such that
$\ell(w_1w_2)=\ell(w_1)+\ell(w_2)$ then we have a canonical
isomorphism \begin{equation}\label{umn} j_{w_1*}*j_{w_2*}\cong
j_{w_1w_2*}. \end{equation}
 If $w_1,w_2,w_3$ are such that
$$\ell(w_1w_2w_3)=\ell(w_1)+\ell(w_2)+\ell(w_3)$$
then the two isomorphisms between $j_{w_1*}*j_{w_2*}*j_{w_3*}$ and
$j_{w_1w_2w_3*}$ arising from \eqref{umn} coincide.

b) $j_{w*}$ is an invertible object of the monoidal category
$D_I$. More precisely, we have $$j_{w*}*j_{w^{-1}!}\cong \delta_e
\cong j_{w^{-1}!}*j_{w*}.\ \ \ \ \ \ \epf$$
\end{Lem}
\begin{Cor} a) The map $\lambda \mapsto j_{\lambda*}$ for
$\lambda\in \Lambda^+$ extends naturally to a monoidal functor
$Rep(T)\to D_I$.

b) The map $V\otimes (\lambda) %\otimes (n)
\mapsto \Z(V)*j_{\lambda*} %(\frac{n}{2})
$ for $V\in Rep(\LG)$,
 $\lambda\in \Lambda^+\subset Rep(\LT)$ %, $n\in \Zet\subset Rep(\Gm)$
extends naturally to a monoidal functor $\Fbar:Rep(\LG\times \LT %\times \Gm
)\to
D_I$.
\end{Cor}
\proofpt The length function on $W$ is additive on the subsemigroup
$\Lambda^+$, thus  Lemma \ref{Lemmult} applies to $w_1,w_2,w_3\in
\Lambda^+$, and implies statement (a). Then (b) follows from the
central property of the functor $\Z$, which yields a commutativity
isomorphism $\Z(V)*j_{\lambda*}\cong j_{\lambda*}*\Z(V)$ satisfying
the pentagon identity. \epf

 We denote the image of $\lambda$ under the functor
 defined in the Corollary by $J_\lambda$. It follows from the
 definition that $J_\lambda=j_{\lambda !}$ for $\lambda\in
 -\Lambda^+$, $J_\lambda = j_{\lambda*}$ for $\lambda \in
 \Lambda^+$, and $J_{\lambda+\mu}\cong J_\lambda*J_\mu$.

Following Mirkovi\' c we call $J_\lambda$  the Wakimoto sheaves.
Theorem \ref{IW} below asserts that $J_\lambda$ are actually {\it
objects
of the abelian category} $\P_I$ (a'priori they are defined as
objects of
the triangulated category $D_I$).

\subsection{Monodromy and "highest weight" arrows: characterization of
the functor $\Ftil$}\label{ar}
\subsubsection{Arrows between perverse sheaves}
 Recall that the monoidal functor $\Z$ comes equipped
with a tensor endomorphism $\M=\{\M_V=\M_{\Z(V)}\in End(\Z(V)) %, \Z(V)(1))
\}$
defined by the logarithm of monodromy (see \cite{KGB}, Theorem 2;
we fix and use an isomorphism $\Ql\cong \Ql(1)$).

We also define an arrow $\b_\lambda:Z_\lambda\to j_{\lambda*}$, $\lambda
\in \Lambda^+$. The definition is clear from the next

\begin{Lem}\label{bFl}
For all $\lambda\in \Lambda^+$ the Schubert cell
$\Fl_{\lambda}$ is open in the support of $Z_\lambda$; and we
have a canonical isomorphism
\begin{equation}\label{jZ}
j_{\lambda}^*(Z_\lambda)\cong
\cons[\ell(\lambda)]. %(\frac{\ell(\lambda)}{2}).
\end{equation}
\end{Lem}
\proofpt Recall that $\pi$ denotes the projection $\Fl\to \Gr$.
It is immediate to see from the definition of the functor
$\Z$ (see
\cite{KGB}, 2.2.3) that the support of $Z_\lambda$ is contained
 in the preimage under $\pi$ of
the closure of the Schubert cell $\Gr_\lambda$; and also that its
dimension equals $\dim \Gr_\lambda=\dim (\Fl_\lambda)$.
Thus it can not contain $\Fl_w$ for $w\succ \lambda$.
 It contains $\Fl_\lambda$, and we have the canonical isomorphism
\eqref{jZ}, because
the support of $\pi_*(Z_\lambda)=IC_\lambda$ contains $\Gr_\lambda$,
and $IC_\lambda|_{\Gr_\lambda}=\cons[\ell(\lambda)] %(\frac{\ell(\lambda)}{2})
$. \epf

\subsubsection{Arrows between coherent sheaves}\label{arr}
First, consider the variety $\Lg$ equipped with the adjoint
action. Then  every $\F\in Coh^\LG(\Lg)$ carries a canonical
endomorphism,
%$N^{taut}$,
 such that the induced endomorphism of the fiber at a
point $x\in \Lg$ coincides with the action of $x\in Stab_{\Lg}(x)$
coming from the equivariant structure; we denote this endomorphism
by $N^{taut}_\F$, and abbreviate $N^{taut}_{V}=N^{taut}_{V\otimes
\O}$.
%We will also consider the category  $Coh^{\LG\times \Gm}(\Lg)$,
%where $\Gm$ acts on $\Lg$ by $t:x\mapsto t^{-2}x$.
%Then we get $N^{taut}_\F\in Hom(\F, \F(2))$, where $\F(n)$ denotes
%the twist of $\F$ by the $n$-th power of the tautological character of $\Gm$.

Next, consider the basic affine space $\LG/\LU$, and its affine
closure $\LGUbar $. We fix an isomorphism between the ring of
regular functions on $\LG/\LU$ and the ring
$\oplusl_{\Lambda_+}V_\lambda$ with multiplication given by
$\oplusl_{\lambda,\mu} m_{\lambda,\mu}$
(see section \ref{mlm} for notation). Then  for $\lambda\in
\Lambda^+$ we get a morphism in $Coh_{fr}^{\LG\times \Gm}(\LGUbar)$:
$$\bo_\lambda:
V_{\lambda} \otimes \O\to \O_\lambda.$$

\begin{Prop}\label{exi}
There exists a unique extension of $\Fbar$ to a monoidal functor
$\Ftil:Coh_{fr}^{\LG\times \LT}(\Nthat_{af})\to \PI$
% (respectively, $\Ftil^{gr}:
%Coh_{fr}^{\LG\times \LT\times \Gm}(\Nthat_{af})\to
%\PI^{gr}$ 
such that
$\Ftil(N_V^{taut}) = \M_{\Z(V)}$, $\Ftil(\bo_\lambda) =
\bFl_\lambda$.
\end{Prop}
The proof of the Proposition will be given at the end of the next
section  after some general nonsense preparation.

\begin{Rem} One can show that any arrow in
 $Coh_{fr}^{\LT}(\Nthat_{af})$ % (or $Coh_{fr}^{\LG\times \LT}(\Nthat_{af})$)
can be obtained from the arrows $B_\lambda$, $N_{V}$ and
identity arrows
by taking tensor products and direct summands. This implies the uniqueness
statement in the Proposition. We will give a slightly different argument
in the next section.
\end{Rem}

The only
geometric statement needed for the proof of Proposition \ref{exi} is the next
\begin{Lem}\label{bFlcommu}
a) For $\lambda, \mu \in \Lambda^+$  we have
$Hom( Z_\mu,J_\lambda)=0$ unless $\lambda\preceq \mu$.

b) For $\lambda,\mu\in \Lambda^+$ the following diagram is
commutative
$$
\CD Z_{\lambda+\mu} @>>{\Z(m_{\lambda,\mu})}> Z_{\lambda}*Z_\mu \\
 @VV{\bFl_{\lambda+\mu}}V  @VV{\bFl_\lambda*\bFl_\mu}V\\
j_{\lambda+\mu*} @= j_{\lambda*}*j_{\mu*}
\endCD
$$
where the lower  horizontal isomorphism  comes from Lemma \ref{Lemmult}(a).

c) We have $\bFl_\lambda\circ \M_{Z_\lambda}=0$.
\end{Lem}
\proofpt a)  As was said in the proof of the previous Lemma,
the support of $Z_\mu$ is contained in the preimage under
$\pi$ of
the closure of the Schubert cell $\Gr_\mu$. It is well known that
$\Fl_\lambda$ is contained in this set iff  $\lambda\preceq \mu$.

b) $\pi_*$ induces an isomorphism of one dimensional vector spaces
$$Hom(Z_\lambda, j_{\lambda*})
\iso Hom(IC_\lambda, \jbar_{\lambda*}),
$$
where $\jbar_{\lambda*}=\jbar_{\lambda*}(
\cons[\ell(\lambda)] %(\frac{\ell(\lambda)}{2}))
$. Thus it suffices
to check that applying $\pi_*$ to the above diagram we get a
commutative one. This follows from the definition, and the
canonical isomorphism  $\pi_*\circ \Z\cong id_{\P_{G(O)}}$.

c) It suffices to see that
$$j_{\lambda}^*(\M_{Z_\lambda})
\in End(j_\lambda^*(Z_\lambda) %,j_\lambda^*(Z_\lambda)(1)
 )=0,$$
 This follows from nilpotency of
$\M_{Z_\lambda}$ and Lemma \ref{bFl} which shows that
$End(j_\lambda^*(Z_\lambda))$ is one dimensional. \epf

\subsection{Tannakian and Drinfeld-Plucker formalism}\label{form}
Notice that the arrows $\bo_\lambda$ introduced in section \ref{arr}
  satisfy the so-called {\it Plucker
relations}, i.e.
\begin{equation}\label{pluck1}
\bo_\lambda\otimes \bo_\mu= \bo_{\lambda+\mu}\circ
(m_{\lambda,\mu}\otimes id_\O).
\end{equation}
%\begin{equation}\label{pluck2}
% (\bo_\lambda\otimes \bo_\mu)\circ (m' \otimes id_\O)
%=0
%\end{equation}
%where
%$m'$ stands for any
% morphism $ V_\nu\to V_\lambda\otimes V_\mu $
% with $\nu\ne \lambda+\mu$.
\begin{Lem}\label{tanak}
 Let $A$ be a commutative algebra with a $\LG$ action.

a) Let $N$ be a tensor endomorphism of the functor $V\mapsto
A\otimes V$; thus $N$ is a collection of  $\LG$-invariant
endomorphisms $N_V\in End_A(A\otimes V)$, $V\in Rep(\LG)$,
functorial in $V$ and such that
$$N_{V_1}\otimes 1 +1 \otimes N_{V_2}=N_{V_1\otimes V_2}$$ for all
$V_1,V_2\in Rep(\LG)$. Then there exists a unique element 
$x_N\in \Lg\otimes A$, such that $N_V$ coincides with the action of $x$
in $A$.

Also, there is a unique $\LG$-equivariant homomorphism 
 $\phi: \O_\Lg\to A$ such that $N_V=\phi_*(N_V^{taut})
\dfn id_A\otimes _{\O_\Lg} N_V^{taut}$.

b) Assume that $A$ is equipped
with a $\Lambda$ grading compatible with the $\LG$ action
(in other words, an action of
$\LT$ commuting with the $\LG$ action is given);
 and suppose that for every
$\lambda \in \Lambda^+$ we are given a $\LG$-equivariant morphism
$\bobo_\lambda: V_{\lambda}\otimes A\to A(\lambda)$ satisfying the
Plucker relations \eqref{pluck1}
(with
$\bo_\lambda$ replaced by $\bobo_\lambda$, and $\O$ replaced by
$A$).
 Then there exists a unique
$\LG\times \LT$-equivariant homomorphism $\phi:
\O({\LG/\LU})=\O(\overline{\LG/\LU})\to A$ such that
$\bobo_\lambda=\phi_* (\bo_\lambda)\dfn  id_A\otimes _{\O_\Lg} \bo_\lambda
$.

c) Let $A, \bobo_\lambda$ be as in (b), and $N$ be as in (a).
Assume that \begin{equation}\label{ravnul}
 \bobo_\lambda \circ  N_{V_\lambda}=0
\end{equation}
for all $\lambda$. Then the homomorphism $\O(\LGU\times \Lg)\to A$
provided by (a,b) factors through $\Ohat_{\Nt}$.
\end{Lem}
\proofpt The first statement in (a) is well-known.
% (cf. e.g.  \cite{Dten}). %??

The second statement in (a) is a restatement of the first one.
More precisely, a homomorphism $\O_{\Lg}\to A$ is specified by an element
of $Hom((\Lg)^*, A)=\Lg\otimes A$, and it is straightforward to see
that $x\in \Lg\otimes A$ satisfies the conditions of the first
statement in (a) iff the corresponding homomorphism 
$\O_\Lg\to A$ satisfies the conditions of the second one.

To check (b) recall that $\O(\LGUbar)\cong \oplusl_{\Lambda^+}
V_\lambda$. Then the requirement on $\phi$ is equivalent
to 
\begin{equation}\label{defphi}
\phi|_{V_\lambda}= b_\lambda|_{V_\lambda\otimes 1}.
\end{equation}
Thus uniqueness of $\phi$ is clear.
The Plucker relations ensure that the map $\phi:\O(\LGU)\to A$
defined by \eqref{defphi} is indeed a homomorphism, which shows existence.

(c) is immediate from the definition of $\Nthat_{af}$.
 \epf

\begin{Prop}\label{FtoC} Let $\C$  be an additive monoidal category.

a) Let  % $H$ be an  algebraic group,  
$F:Rep(\LG)\to \C$ be a monoidal functor.

 Let $N=\{N_V\}$ be a tensor endomorphism of $F$, such that the
image under $F$ of the commutativity isomorphism in $Rep(\LG)$  is
functorial with respect to $N$. % Assume that $H$ is reductive.
 Then
there exists a unique extension of $F$ to a monoidal functor
$\Ftil: Coh_{fr}^{\LG}(\Lg)\to \C$ such that $N_V={\tilde F}(N_V^{taut})$ for
all $V\in Rep(\LG)$ (here $\LG$ acts on $\Lg$ by the
adjoint action).

b) Let $F:\LG\times \LT\to \C$  be a monoidal functor.
 Suppose that for each  $\lambda\in \Lambda^+$ 
we are given
 transformations $\bobobo_\lambda:F(\lambda)
\to F(V_\lambda)$ satisfying the
Plucker relations, i.e. such that
$$ %\begin{equation}\label{pluck1}
\bobobo_\lambda\otimes \bobobo_\mu= \bobobo_{\lambda+\mu}\circ
F(m_{\lambda,\mu}).
 $$ % \end{equation}

%$$ %\begin{equation}\label{pluck2}
% (\bobobo_\lambda\otimes \bobobo_\mu)
%\circ F(m')=0.
%$$
%where $m_{\lambda,\mu}$  are as in Lemma \ref{tanak}.

 Assume
that  the image of the commutativity isomorphism under $F$ is
functorial with respect to $\bobobo_\lambda$.
 Then
there exists a unique extension of $F$ to a monoidal functor
$\tilde F:Coh_{fr}^{\LG\times\LT}(\overline{\LG/\LU})\to \C$ such that
$\bobobo_\lambda=\tilde F(\bo_\lambda)$.

c) Let $F$, $\bobobo_\lambda$ be as in (b), and $N\in
End(F|_{Rep(\LG)})$ be as in (a).
Assume  that
\begin{equation}\label{ravnulka}
\bobobo_\lambda \circ N_{V_\lambda} = 0
\end{equation}
for all $\lambda$. Then (a,b) provide an extension of $F$ to a
monoidal functor $Coh_{fr}^{\LG\times\LT}(\LG/\LU\times \Lg)\to \C$,
 which factors
through $\Ohat_{\Nt_G}-mod _{fr}^{\LG\times \LT}$.
\end{Prop}

\proofpt Let $H$ stand for $\LG$ if we are in the situation of (a),
and for $\LG\times \LT$ if we are in the situation of either
(b) or (c).

 First we claim that without loss of generality we can
assume that $\C$ is a tensor category, and $F$ is a tensor
functor. More precisely, we claim that
% in the situation of point (a) %(respectively, (b) or (c))
 it is possible to factor $F$ as a
composition $F=F''\circ F'$, where $F'$ is a tensor functor from
$Rep(H)$ to a tensor category $\C'$, and $F'':\C'\to \C$ is a
monoidal functor; moreover, in the situation of (a)
there exists a tensor endomorphism
$N'$ of $F'$
%(respectively, morphisms $\bobobo_\lambda':F'(k_\lambda)\to
%F'(V_\lambda) $, or both $N'$ and $\bobobo_\lambda'$)
satisfying the conditions of (a), % (respectively, (b) or (c))
 and
such that $N=F''(N')$; and similarly for (b) and (c).
 %(respectively, $\bobobo_\lambda=F''(\bobobo_\lambda')$).

Namely, we can define $\C'$ as follows. We set
$Ob(\C')=Ob(Rep(H))$, and $Hom_{\C'}(V,W)\subset Hom_\C(F(V),F(W))
$ consists of such elements $\phi$ that for all $U\in Rep(H)$ the
diagram
$$\CD
F(U)\otimes F(V)  @>{id\otimes \phi}>> F(U)\otimes F(W) \\
@VVV  @VVV \\
F(V)\otimes F(U) @>{\phi\otimes id}>> F(W)\otimes F(U)
\endCD
$$
is commutative for all $U\in Rep(H)$; here the vertical arrows are
images under $F$ of the commutativity isomorphism. The pentagon
identity implies that $\C'$ is indeed a tensor category; the
definition of $F'$, $F''$, $N'$, $\bobobo_\lambda'$ is clear.
So from now on we will assume that $\C$, $F$ are tensor.

We will use underlined symbols to denote representations of $H$, and the corresponding
unadorned symbol will denote the underlying vector space.

 Let   $\unO$ be the  module of regular functions
on $H$ where $H$ acts by left translations; thus $\unO$ is an
ind-object of $Rep(H)$, $\unO=\oplusl_{\unV\in IrrRep(H)}
V^*\otimes \unV $, where $IrrRep(H)$ is a set of representatives
for isomorphism classes of irreducible representations of $H$.
Thus $\unO$ is a commutative ring ind-object in $Rep(H)$.
Set $$A=Hom_\C(1_\C,F(\unO))=\oplusl_{\unV\in IrrRep(H)}
V^*\otimes
Hom_\C(1_\C,F(\unV));$$ here (and below) we use the same notation
for a functor on a category, and the induced functor on the
category of ind-objects.
Then $A$ is an associative algebra equipped with an $H$ action.
Commutativity of $\unO$ and tensor property of $F$ show that $A$
is commutative; and the functor $\Phi:X\mapsto Hom_\C(1, X\otimes
F(\unO))$ is a tensor functor from the full image of $F$ to
$A-mod_{fr}^H$. It is easy to see from the definitions that $\Phi$
induces an isomorphism $$Hom_\C(1_\C, F(\unO))\iso
A=Hom_{A-mod_{fr}^H}\left(\Phi(1_\C),
\Phi\left(F(\unO)\right)\right).$$ Since $H$ is reductive,
$Rep(H)$
is semisimple, and every irreducible $\unV\in Rep(H)$ is a direct
summand of $\unO$; hence for all $\unV$ we have an isomorphism induced
by $\Phi$
$$Hom_\C(1_\C, F(\unV))\iso  Hom _{A-mod^H_{fr}}(\Phi(1_\C), \Phi(F(\unV)) ) \cong
Hom _{A-mod^H_{fr}}(A, V\otimes A).$$
Since $Rep(H)$ is rigid, we see that $\Phi$ is a full embedding.
Thus we can assume that $\C=A-mod^H$ for a commutative algebra with an
$H$-action, and $F:\unV \mapsto V\otimes A$. In this case
the statements of the Proposition reduce to that of Lemma
\ref{tanak}. \epf

\subsubsection{Proof of Proposition \ref{exi} and definition of
the functor $\Ftil$}
The Proposition follows from Lemma \ref{bFlcommu} in view
of Proposition \ref{FtoC}.

\medskip

We now extend $\Ftil$ to the homotopy category
 $Hot(Coh_{fr}^{\LG\times \LT}(\Nthat_{af}))$.

\subsection{Proof of Proposition \ref{Frq}}
We will construct an isomorphism
$$Fr\circ \Ftil\iso \Ftil\circ \bq^*;$$
the claim about $F$ follows (once we show
that $F$ exists).

Uniqueness part of Proposition \ref{FtoC} shows that we will be done if we 
construct an isomorphism of monoidal functors 
%\begin{equation}\label{onRep}
$$
\phi:Fr\circ \Fbar\iso \Fbar
$$
%\end{equation}
such that
\begin{equation}\label{onRep} 
\begin{array}{ll}
\phi(Fr^*(\bo_\lambda))=\bo_\lambda, \\
\phi(Fr^*(\M_Z(V)))=q^{-1}\cdot \M_{\Z(V)}. 
\end{array}
\end{equation}

An isomorphism $\phi$ induces a structure of a Weil sheaf on
 $J_\lambda$, $Z_\mu$, and it is clearly uniquely determined
by this structure. For $\lambda\in \Lambda^+$ we fix the Weil structure
on $J_\lambda=j_{\lambda*}$ so that the resulting Weil sheaf is 
$j_{\lambda*}\cons[\ell(w)](\lT)$. We also require that the isomorphism
$J_\lambda*J_\mu\cong J_{\lambda+\mu}$ lifts to an isomorphism of Weil
sheaves; this fixes the Weil structure on $J_\lambda$ for all $\lambda$.

Let us now define the Weil sheaf which provides the desired
isomorphism $Fr^*(Z_\lambda)\cong Z_\lambda$.
The functor $\Z:\P_\K(\Gr)\to \P_I(\Fl)$ is actually defined as a functor
between the categories of Weil sheaves (if one fixes the splitting
of the surjection $Gal(\Fq((t)))\to Gal(\Fq)$, cf. the footnote on p. 263
in \cite{KGB}). Then Weil sheaf in question is defined to be 
$Z_\lambda^{Weil}:=\Z(IC_\lambda^{Weil})$
where $IC_\lambda^{Weil}=\jbar_{\lambda !*}\left(
\cons[\ell(\lambda)]        (\frac{\ell(\lambda)}{2})
\right) $. 

These requirements clearly define the tensor isomorphism $\phi$ uniquely.
Verification of existence of $\phi$ reduces to checking
 that the
isomorphism 
\begin{equation}\label{SaF}
Fr\circ \Sa\cong \Sa,
\end{equation} providing $IC_\lambda$ with the
Weil structure isomorphic to $IC_\lambda^{Weil}$ is tensor; the rest
then follows from $\Z$ being tensor. Existence of a tensor structure on
\eqref{SaF} will be clear if we show that the convolution $IC_\lambda^{Weil}
*IC_\mu^{Weil}$ is isomorphic to a direct sum of Weil sheaves
$IC_\nu^{Weil}$. We now prove this.

 Notice that 
 Frobenius acts on 
the  total cohomology $H^\bu(IC_\lambda^{Weil})$
 by a diagonalizable automorphism
with eigenvalues $q^{n/2}$, $n\in \Zet$; this follows e.g. from \cite{BGS},
\S 4.4. The 
 functor of   total cohomology
on $\P_\K(\Gr)$
carries a tensor structure (see \cite{MV}, \cite{BD}); %??
the latter is readily seen to be compatible with the Frobenius action. 
Thus the action of Frobenius on $H^\bu(IC_\lambda^{Weil}*IC_\mu^{Weil})$
is diagonalizable with  eigenvalues $q^{n/2}$. This implies the desired
statement, because we know that $IC_\lambda*IC_\mu\cong \oplus IC_\nu$,
 and the action of Frobenius
 on cohomology determines the isomorphism
class of a Weil sheaf which is geometrically isomorphic to a direct sum
of $IC_\nu$, $\nu\in \Lambda^+$.

It remains to check \eqref{onRep}.
The first equality in \eqref{onRep} is clear from the definition.
The second one follows from the fact that for
an $l$-adic sheaf $\F$ the logarithm of monodromy 
on nearby cycles is a morphism of Weil sheaves
$\Psi(\F)\to \Psi(\F)(-1)$. \epf

\subsection{Filtration of central sheaves by Wakimoto sheaves}
\label{filtrWak}
%{\it In this section we work over algebraically closed field
%$k$.} 

 The property of central sheaves proved in this
section
is a geometric counterpart of Bernstein's description
of the center $Z_H$ of the Iwahori-Matsumoto Hecke algebra $\H$, which says
that $Z_H=\Ce[\theta_\lambda]^{W_f}$; moreover, the map
$K(Rep(\LG))\to Z_H$ sends the class of representation $V$ to its
character $\chi_V\in k[\Lambda]=k[\theta_\lambda]$ (see e.g. \cite{doistor},
Theorem 8.1).
Bernstein presentation for $H$ (in particular, the elements
$\theta_\lambda$) can be easily described in terms of their action
in the space of $I$-invariant vectors in the universal principal
series representation $C_c(G(F)/(N(F)\cdot T(O))$, and thus in
terms of their integrals over $N(F)$-orbits in $G(F)$. Quite
similarly, the property of central sheaves proved in this section
is related to computation of compactly supported cohomology of
their restrictions to $\NF$-orbits.

The next Theorem, which is the main result of this section,
contains two close statements. Statement (a) will be used later;
statement (b) is included for completeness.

 Recall that the $\NF$-orbits on $\Fl$ are parameterized by $W$;
and $\io_w:\Si_w\imbed \Fl$ denotes the embedding of an orbit.

We fix a total ordering $\leq$ on the group $\Lambda$ compatible with the
group structure and with the standard partial ordering (i.e.
$\lambda >\mu$ if $\lambda-\mu$ is a sum of positive roots).
\begin{Thm}\label{filtrWakr}
a) For $V\in Rep(\LG)$ the sheaf $ \Z(V)$ has a unique filtration
indexed by $(\Lambda, \leq)$ such that the associated graded
$gr_\nu(\Z(V))
=\Z(V)_{\leq \nu}/\Z(V)_{<\nu}$
is of the
form
$$gr_\nu(\Z(V))\cong J_\nu\otimes W_V^\nu$$ for some vector space
$W_V^\nu$. The functor $$\Phi:V\mapsto \oplusl _\nu W_V^\nu=
\oplusl _\nu Hom(J_\nu, gr_\nu(\Z(V) )$$ is a
tensor functor from $Rep(\LG)$ to the category of $\Lambda$-graded
vector spaces (obviously equivalent to $Rep(\LT)$).
$\Phi$ is isomorphic the restriction functor $Rep(\LG)\to
Rep(\LT)$;
 in particular, $\dim
W_V^\nu$ equals the multiplicity of the weight $\nu$ in $V$.

b) The space   $ H^i_c(\io_{w}^*(\Z(V)))$ vanishes
 unless $w=\nu\in \Lambda$,
 $i= \ell(\nu)$; in which case we have
$$H^{\ell(\nu)}(\Z(V))\cong W_V^\nu$$
where $W_V^\nu$ is as in (a).
\end{Thm}
\begin{Rem} The sheaf $V_\mu \otimes \O_{\Nt}$ carries a
filtration with subquotients being sums of line bundles (this
filtration is actually a pull-back of a filtration on $V_\mu
\otimes \O_{\LG/\LB}$). It will be clear from the
construction of the functor $F$ that  the filtration of Theorem
\ref{filtrWakr}(a) is the image of this filtration under  $F$.
\end{Rem}
\begin{Rem} In \cite{MV} Mirkovi\' c and Vilonen prove a result
similar to part (b) of the above Theorem; namely, they compute the
compactly supported cohomology of $N(F)$ orbits with coefficient
in an irreducible object of $\P_\K(\Gr)$. One can show that the
two results are actually equivalent.
\end{Rem}
The proof of the Theorem occupies the rest of this section.
\subsubsection{(a) implies (b)}\label{gammatam}
The last statement in the next Lemma yields the
implication (a) $\Rightarrow$ (b).

\begin{Lem}\label{sdvig}
%a)
 For $\lambda\in \Lambda$, and $X\in D_I$ we have
 \begin{equation}\label{sdvigeq}
H^\bu(\io_{\lambda\cdot w}^!J_\lambda*X)=
H^\bu( \io_w^!(X))[(\lambda,2\rho)].
\end{equation}
In particular, $H^i( \io_w^!(J_\lambda))=0$ unless $w=\lambda$,
 $i=\ell(\lambda)$, in which case it has dimension one.
%b)  For $w\in \Lambda^+\cdot W_f$ (in particular, for $w=\lambda\in \Lambda^+$)
%we have $\Fl_w\subset \Si_w$.
\end{Lem}
\proofpt 
It is clear that if \eqref{sdvigeq} holds for $\lambda_1$,
$\lambda_2$
then it also holds for $\lambda_1-\lambda_2$. Thus
we can assume without loss
of generality that $\lambda \in \Lambda^+$.

 For $w\in W$ let $\wtil$ be a
representative of the coset $w\in Norm(T(O))/T(O)$,
where $Norm(T(O))$ is the normalizer of $T(O)$. It follows from the
definitions that for $X\in D_I$ we have $$j_{w*}*X\cong \Gamma
_{\bI\cap \wtil \bI \wtil^{-1}} ^\bI \wtil_{*}(X)[\ell(w)].$$
(Notations for the induction functor $\Gamma$ were recalled
before  before Lemma \ref{Hpo} above.) It
is
clear that for $\lambda\in \Lambda$
$$H^\bu(\io_{\lambda\cdot w}^!\tilde\lambda_*(X))= H^\bu(
\io_w^!(X))$$ since
$\tilde\lambda(\Si_w)=\Si_{\lambda\cdot w}$. Also it is not
difficult to check that for %$w\in \Lambda^+\cdot \Wf$, in particular for
$\lambda\in \Lambda^+$ we have $\lambdatil \bI \lambdatil^{-1}
\supset
\bI\cap \BFm$. Then the triangular decomposition $\bI=\bI\cap\NF
\cdot \bI\cap \BFm$ yields an isomorphism $$ \Gamma_{\bI\cap
\lambdatil
\bI \lambdatil^{-1}}^\bI = \Gamma_{\bI\cap \lambdatil \bI
\lambdatil^{-1}\cap \NF
}^{\bI\cap \NF}.$$ The induction functor $\Gamma_{H'}^H$ commutes
with the $!$ restriction to an $H$-invariant subvariety; when $H$,
$H'$ are unipotent it also does not change the total cohomology.
Applying this
observation to $H=\bI\cap \NF$, $H'=\bI\cap \lambdatil \bI
\lambdatil^{-1}\cap \NF$ and the subvariety $\Si_{\lambda\cdot
w}\subset \Fl$ we get the statement. \epf
\subsubsection{Uniqueness of the filtration}\label{unifi}
 Uniqueness of the
filtration follows from the following
\begin{Lem}
We have $Hom^\bu(J_\lambda, J_\mu)=0$ unless $\lambda \preceq
\mu$; and $Hom^\bu(J_\lambda,J_\lambda)=\Ql$.
\end{Lem}
\proofpt  Pick $\nu$ such that $\nu+\lambda, \nu+\mu \in
\Lambda^+$. Since the functor of convolution with $J_\nu$ is
invertible we have
$$Hom^\bu(J_\lambda, J_\mu)=Hom^\bu(J_{\nu+\lambda}, J_{\nu+\mu})=
Hom^\bu(j_{\nu+\lambda*}, j_{\nu+\mu*}).$$ The latter space can be
nonzero only if $\Fl_{\nu+\lambda}$ lies in the closure of
$\Fl_{\nu+\mu}$, which is known to be equivalent to
$\lambda\preceq \mu$. \epf
\subsubsection{Existence of the filtration}
We will say that an object
 $X\in \P_I$ is {\it convolution exact} if
 $X*L\in \P_I$ for all $L\in \P_I$.
We will say that $X$ is {\it central} if $X*L\cong L*X$ for all
$X\in \P_I$.

It will be convenient to extend the definition of $J_\lambda$ to
all $w\in W$
by setting $J_w=J_\lambda* j_{w_f*}$ for $w=\lambda\cdot w_f$,
$\lambda\in
\Lambda$, $w_f\in W_f$.

The next  result  is proved in the Appendix.

 \begin{Thm}\label{IW} a) The objects $J_w\in D_I$ actually lie in $\P_I$.

b) $\Fl_w$ is open in the support of $J_w$, and $j_w^*(J_w)\cong
\cons [\ell(w)]$.
 \end{Thm}

 The next Proposition obviously implies the existence
of the filtration.
\begin{Prop}\label{fi}
a) Any convolution exact object of $\P_I$ has a filtration whose
subquotients
are Wakimoto sheaves $J_w$.

b) If $X$ is also central then only $J_w$ with $w\in \Lambda$
appear in the
filtration of (a).
\end{Prop}
\begin{Rem}\label{remtilt}
Statement (a) of the Theorem can be compared to the following result
due (to the best of our knowledge) to Mirkovi\' c (unpublished): 
every convolution exact sheaf on the finite
dimensional flag variety $G/B$ which is smooth along the Schubert
stratification is {\it tilting}, i.e. has a filtration with 
subquotients $j_{w!}$, and also a  filtration with
subquotients $j_{w*}$.

We sketch a proof of Mirkovi\' c's result for the sake of completeness.
Let $\F$ be a convolution exact perverse sheaf on $G/B$ as above. We have to check
that $Ext^{>0}(j_{w!}, \F)=0=Ext^{>0}(\F,j_{w*})$. We check the first equality, the other one is similar.
Since $\F$ is convolution exact, the convolution $\F*j_{w_0!}$ is a perverse sheaf, thus it lies
in the full subcategory generated by the objects $j_{w!}[d]$, $w\in W$, $d\geq 0$ under extensions.
Thus $\F=\F*j_{w_0!}*j_{w_0*}$ lies in the full subcategory generated by
$j_{w!}*j_{w_0*}[d]=j_{ww_0*}[d]$, $d\geq 0$, which implies the needed Ext vanishing.

The central sheaves $Z_\lambda$ (for
$\lambda\ne 0$) provide
examples of convolution exact objects of $\P_I$ which are not tilting (see, however,
Theorem \ref{tilt} and Remark \ref{remtiltasph} below).
\end{Rem}
The proof of the Proposition will be given after some auxiliary
Lemmas.
\begin{Lem}\label{onWak}
a) We have $J_\lambda*J_w\cong J_{\lambda \cdot w}$.

b) If $w\in \Lambda^+\cdot W_f$ then $J_w= j_{w*}$.
If $w\in (-\Lambda^{++})\cdot W_f$ then $J_w=j_{w!}$; here
$\Lambda^{++}$
is the set of strictly dominant weights.
\end{Lem}
\proofpt (a) is immediate from the definitions. To prove (b) we
observe
that for $w_f\in W_f$, $\lambda\in
\Lambda^+$, $\mu\in  \Lambda^{++}$ we have
$$\ell(\lambda\cdot w_f)=\ell(\lambda)+\ell (w_f) \Rightarrow
j_{\lambda\cdot w_f*}=j_{\lambda*}*j_{w_f*}=J_\lambda* J_{w_f}=
J_{\lambda\cdot w_f};$$
$$\ell((-\mu)\cdot w_f)=\ell(-\mu)-\ell (w_f) \Rightarrow
j_{-\mu\cdot w_f!}=j_{-\mu!}*j_{w_f*}=J_{-\mu}* J_{w_f}=
J_{-\mu\cdot w_f}. \ \ \epf$$

\subsubsection{Perverse sheaves on stratified spaces}\label{genizvr}
We now recall some facts about perverse sheaves on stratified
spaces.

 Let $X=\cupl_{s\in S} X_s$ be a
stratified scheme over a field; thus $X_s\subset X$ are locally closed smooth
subschemes. 
We assume for simplicity of notations
 that the embeddings $j_s:X_s\imbed X$ are affine, and that
$j_{u}^*j_{s*}(\cons)$ has constant cohomology sheaves for all $u,s\in S$. We
abbreviate $j_{s*}=j_{s*}(\cons[\dim X_s])$,
$j_{s!}=j_{s!}(\cons[\dim X_s])$. Let $D$ be the derived category
of constructible sheaves on $X$, and $(D^{<0},\,D^{>0} )$ be the
perverse $t$-structure, and $\P$ be its heart
 (the category of perverse sheaves).
The following statement is standard.
\begin{Claim}\label{obsh}
For $\F\in D$ set $S_\F^*=\{s\in S \ |\ j_s^*(\F)\ne 0\}$;
 $S_\F^!=\{s\in S \ |\ j_s^!(\F)\ne 0\}$.
We have

a) If $\F\in D^{\leq 0}$ and cohomology of $j_s^*(X)$ are constant
sheaves for all $s\in S$, then $\F\in \< j_{s!}[i]\ | i\geq 0,
s\in S_\F^* \rangle$ (cf section \ref{prooW} for notations).

b) If $\F\in D^{\geq 0}$ and $\F\in \< j_{s*}[i]\ | i\geq 0, s\in
S \rangle$ then $\F$ is a perverse sheaf; moreover, $\F$
carries a filtration with subquotients
isomorphic to $j_{s*}$, $s\in S_\F^!$.
\end{Claim}
\proofpt (a) is equivalent to saying that $H^k(j_s^*(\F))=0$ for
$k>-\dim (X_s)$,
and is constant otherwise. Here the first condition is the
definition
of the perverse $t$-structure, and the second one was imposed as
an assumption.

The assumptions of (b) imply that $H^k(j_s^!(\F))=0$ for $k>-\dim
(X_s)$,
and is constant otherwise. However, $H^k(j_s^!(\F))=0$
 for $k<-\dim (X_s)$ by the definition of the perverse $t$-structure.
Hence $j_s^!(\F)\cong \cons[\dim (X_s)]^{\oplus n}$, which implies
the
conclusion of (b). \epf

For $X\in D_I$ set
$$W_{X}^*=\{w\in W\ |\ j_w^*(X)\ne 0\};$$
 $$W_{X}^!=\{w\in W\ |\ j_w^!(X)\ne 0\}.$$
\begin{Lem}\label{S} For $X\in D_I$
 there exists a finite subset $S\subset W$, such that
for all $w\in W$ we have
$$W_{j_{w*}*X}^!,\; W_{j_{w!}*X}^* \;\subset w\cdot S;$$
$$W_{X*j_{w*}}^!,\; W_{X*j_{w!}}^*\subset S\cdot w.$$
\end{Lem}
\proofpt
Proper base change shows that any point $x\in \Fl$ such that the
stalk of
 $j_{w!}*X$ at $x$ is nonzero lies in
 the convolution of sets $\Fl_w$ and $Supp(X)$
(i.e. in the image of $\Fl_w\underset{\bI}{\times}Supp(X)$ under
the convolution map). 
Thus to prove the first of the four statements it is enough to show 
the corresponding estimate for convolution of sets; the other three statements follow in a similar way.
Thus for a fixed
$\bI$-invariant ${\mathfrak S}\subset \Fl$
 we have to show that for some $S\subset W$, the convolution of sets 
 ${\mathfrak S}*\Fl_w$ (respectively, $\Fl_w* {\mathfrak S}$) is contained in 
 $\cupl_{v\in S\cdot w} \Fl_w$ (respectively, $\cupl_{v\in w\cdot S} \Fl_w$). Without loss of generality
 we can assume that ${\mathfrak S} = \Fl_v$ for some $v\in W$. The claim easily follows by induction 
 in $\ell(v)$. \epf
 
\proof of Proposition \ref{fi}. a) Let $X\in \P_I$ be convolution
exact, and let
 $S$ be as in  Lemma \ref{S}.
 We can write $S$ as
$$S=\{\lambda_i w_i\},$$ $w_i\in W_f$, $\lambda_i\in
\Lambda$. Choose $\nu_0\in -\Lambda^+$
 such that $\nu_0+\lambda_i\in -\Lambda^{++}$;
thus $j_{(\nu_0\cdot s)!}=J_{(\nu_0\cdot s)!}$ for $s\in S$. Since
$j_{-\nu_0!}*X=J_{-\nu_0}*X\in \P$
we see by Lemma \ref{S}, Claim \ref{obsh}(a) that
$$J_{-\nu_0}*X\in \< j_{(-\nu_0)\cdot s!}[i]\ |\ i\geq 0, s\in S\rangle=
 \< J_{(-\nu_0)\cdot s!}[i]\ |\ i\geq 0, s\in S\rangle.$$
Hence $$J_\nu*X\in \< J_{\nu\cdot s!}[i]\ |\ i\geq 0, s\in
S\rangle$$ for all $\nu\in \Lambda$. In particular, choosing
$\nu\in \Lambda^+$ such that $\nu+\lambda_i\in \Lambda^+$ we see
that $$J_\nu* X\in \< j_{\nu \cdot s*}[i]\ |\ i\geq 0, s\in
S\rangle.$$ Since $J_\nu*X\in \P_I$ this implies statement (a) by
Claim \ref{obsh}(b).

(b) We can choose $\lambda\in \Lambda^+$ such that $\lambda\cdot S
\subset \Lambda^{++}\cdot W_f$, $S\cdot \lambda\subset W_f\cdot
\Lambda^+$. Then we see that $$W_{X*J_\lambda}^!\subset W_f\cdot
\Lambda^{++}\cap \Lambda^{++}\cdot W_f=\Lambda^{++},$$ which
implies statement (b). \epf
\subsubsection{Construction of tensor structure on $\Phi$}\label{A}
Let $\A \subset \P_I$ be the full subcategory of sheaves which
admit a filtration whose subquotients are Wakimoto sheaves
$J_\lambda$
(which makes sense by Theorem \ref{IW}(a)). Since $J_\lambda*J_\mu
=J_{\lambda+\mu}$ we see that $\A$ is a monoidal subcategory of $D_I$.

Let $\grA \subset \A$ be the subcategory whose
objects are sums of sheaves $J_\lambda$, and morphisms are direct
sums of isomorphisms $J_\lambda\to J_\lambda$ and zero arrows.
Thus $\A$, $\grA$ are monoidal subcategories in $D_I$, and $\grA$
is obviously equivalent to $Rep(\LT)$. Since
$Ext^1(J_\lambda,J_\mu)=0$ for $\mu \not \prec\lambda$ (in
particular, for $\mu \geq \lambda$)
 every object  $X\in\A$
actually admits a filtration $(X_{\leq \nu})$ indexed by
$(\Lambda, <)$ such that $gr_\nu (X)=X_{\leq \nu}/X_{<\nu}$ is the
sum of several copies of $J_\nu$. (Recall that $\geq$ is some```
complete order on $\Lambda$ compatible with the standard partial order).
 Since $Hom(J_\lambda,J_\mu)=0$
for $\lambda \not \preceq \mu$, in particular, for $\mu>\lambda$,
such filtration is unique. Thus taking the associated graded is a
well defined functor $gr:\A \to \grA$.

The next statement is an equivalent form of Theorem
\ref{filtrWakr}(a).

\begin{Thm}\label{grfib}
The functor $gr\circ \Z:Rep(\LG)\to \grA\cong Rep(\LT)$ is tensor,
and is isomorphic to the functor of restriction to a maximal
torus.
\end{Thm}

The proof of the Theorem will be given at the end of the subsection.

\begin{Prop}\label{moncen}
a) The functor $gr:\A \to \grA$ has a natural monoidal structure.

b) The composition $gr\circ \Z:Rep(\LG)\to \grA$ has a natural
structure of a central functor (see  \ref{mlm} for
the definition of a central functor).
\end{Prop}

\begin{Lem}\label{fiLem}
Let $D, \otimes$ be a triangulated monoidal category
(where $\otimes$ is triangulated in each variable), and
$\A\subset D$ be a heart of a $t$-structure. Let $A,B\in \A$ be
objects with  filtrations $(A_{\leq i}, B_{\leq i})$. Assume that
$gr(A)\otimes gr(B)\in \A$. Then $A_{\leq i}\otimes B_{\leq j}\in
\A $; and we have a natural isomorphism \begin{equation}
\label{natism}gr(A\otimes B) \cong gr(A)\otimes gr(B),
\end{equation}
where $gr(A\otimes B)$ is the associated graded with respect to
the tensor product filtration $(A\otimes B)_{\leq k}=\suml_i
A_{\leq i}\otimes B_{\leq k-i}$. For a third filtered object $C\in
\A$ the isomorphism \eqref{natism}  is compatible with the
associativity isomorphism.
\end{Lem}

\proofpt The first statement is obvious. To see the second one
notice that for all $i,j$ the morphism $A_{\leq i}\otimes B_{\leq
j}\to (A\otimes B)_{\leq i+j}$ factors through  an arrow
$s_{i,j}:gr_i(A)\otimes gr_j(B)\to gr_{i+j}(A\otimes B)$. Also the
image of $(A\otimes B)_{\leq i+j}$ in $(A/(A_{<i}))\otimes
(B/B_{<j})$ equals $gr(A)_i\otimes gr(B)_j$ which induces an arrow
$\sigma_{i,j}: gr(A\otimes B)_{i+j}\to gr(A)_i\otimes gr(B)_j$. It
is clear that $\sigma_{i,j}\circ s_{i,j}=id$, and that
$$s=\suml_{i,j} s_{i,j}:gr(A)\otimes gr(B)\to gr(A\otimes B)$$
is surjective. Hence $s$ is an isomorphism. Compatibility with
associativity is clear. \epf

\begin{Lem}\label{cencomp} Let $F:\cT\to \C$ be a central functor from a tensor
category $\cT$ to a monoidal category $\C$. Let $G:\C\to \C'$
be a monoidal functor to another monoidal category $\C'$. Assume
that
$G$ admits a right inverse, i.e.
there exists a monoidal functor $G':\C'\to \C$ such that $G\circ
G'\cong id$.
Then $G\circ F$ is naturally a central functor.
\end{Lem}

\proofpt Let $\sigma_{X,Y}:F(X)\otimes Y\to Y\otimes F(X)$,  $X\in
\cT$, $Y\in \C$ be the centrality isomorphism for $F$. Define the
centrality isomorphism for $G\circ F$ by
$\sigma'_{X,Y}=G(\sigma_{X,G'(Y)})$, $X\in \cT$, $Y\in \C'$. Then
$\sigma'$ provides $G\circ F$ with a structure of a central
functor, because all the required diagrams commute being images of
commutative diagrams in $\C$. \epf

\proof of Proposition \ref{moncen}. (a) is immediate from Lemma
\ref{fiLem}.  (b) follows from Lemma \ref{cencomp} by
setting $\cT=Rep(\LG)$, $\C=\A$, $\C'=\grA$, $G=gr$, $G'$ is the embedding
$\grA\imbed \A$. \epf

To prove Theorem \ref{grfib} we need another
\begin{Lem}\label{everycen}
Let ${\cT_1}$, ${\cT_2}$ be abelian rigid tensor categories, and $F:{\cT_1}\to {\cT_2}$,
 $\{\sigma_{X,Y}
:F(X)\otimes Y\iso Y\otimes F(X)\}$ be an additive % {\bf is it automatic?}
 central functor.  Let  $U\in {\cT_2}$ be an object. Suppose that there exist
$V\in {\cT_1}$, and a surjective
map $f:F(V)\to U$ such that for all $X\in {\cT_2}$ the following diagram
is commutative
$$
\CD
F(V)\otimes X @>{\sigma_{V,X}}>> X\otimes F(V)\\
@VVV @VVV \\
U\otimes X @>{C_{U,X}^{\cT_2}}>> X\otimes U,
\endCD
$$
where $C^{\cT_2}$ denotes the commutativity isomorphism in ${\cT_2}$.
Then $\sigma_{X,U}=C_{F(X),U}$ for all $X\in {\cT_1}$.
\end{Lem}
\proofpt By the definition of a central functor we have
%\begin{equation}\label{dia}
$$\sigma_{V,F(X)}\circ \sigma_{X,F(V)}=F(C^{\cT_1}_{V,X}\circ
C^{\cT_1}_{X,V})=id_{F(V)\otimes F(X)},
$$ %\end{equation}
where $C^{\cT_1}$ is the commutativity isomorphism in ${\cT_1}$. The
morphism $C_{U,X}\circ \sigma_{X,U}$ is a quotient of
 $\sigma_{V,F(X)}\circ \sigma_{X,F(V)}$ (recall that tensor product
in a rigid abelian tensor category is exact in each variable);
hence  $C_{U,X}^{\cT_2}\circ \sigma_{X,U}=id_{F(X)\otimes U}$, and
$ \sigma_{X,U}=C_{X,U}^{\cT_2}$. \epf

\subsubsection{Proof of Theorem \ref{grfib}}
 By Proposition \ref{moncen}
 the functor $gr\circ \Z:Rep(\LG)\to \grA\cong
Rep(\LT)$ is central. We need to check that it is in fact tensor.
It suffices to check
 that
\begin{equation}\label{r}
\sigma_{V,A}=C_{gr(\Z(V)),A}^{gr(\A)}
\end{equation}
for $V\in Rep(\LG)$, $A\in gr(\A)$; the tensor property
would then follow from the definition of a central functor.
 Lemma \ref{bFlcommu}(b)
implies that conditions of Lemma  \ref{everycen} hold for $\cT_1=Rep(\LG)$,
$\cT_2=gr\A$,
$V=V_\lambda$, $U=J_\lambda$, $\lambda\in \Lambda^+$.
 Hence \eqref{r} holds for $A=J_\lambda$, $\lambda\in \Lambda^+$.
If \eqref{r} holds for some (rigid) object $A$ then its validity for
another object $A'$ is equivalent to its validity for $A\otimes A'$.
Thus \eqref{r} holds for all $V$, $J_\lambda$, $\lambda\in
\Lambda$; and hence holds always.

% It is clear that $gr\circ \Z$ is faithful, hence
Thus  $gr\circ \Z$  comes from
a homomorphism of algebraic groups $\LT\to \LG$. This homomorphism
is injective, because for every $\lambda\in \Lambda^+$ the
character $\lambda$ is a direct summand in $gr\circ
\Z(V_\lambda)$ by Lemma \ref{bFl}. Hence the image of $\LT$ in
$\LG$ under the above homomorphism is indeed a maximal torus. 

This establishes Theorem \ref{grfib} and thus also Theorem \ref{filtrWakr}.
\epf

\subsection{Factoring $\Ftil$ to  $F$}\label{fa}
Let $\de \Nthat\subset \Nthat_{af}$ be the complement to $\Nthat$.
We will show that $\Ftil %, \Ftil^{gr}
$ yields a functor
$D^\LG(\Nt)\to D(\A)$
%, $D^{\LG\times \Gm}(\Nt)\to D(\A^{gr})$
by checking that it sends all complexes whose cohomology is
supported on
$\de \Nthat$ to acyclic complexes.
 This will be deduced from the existence of a
filtration
on $Z_\lambda$ constructed in the previous section (recall that
the definition
of $\Ftil$ only relied on Lemmas \ref{bFl}, \ref{bFlcommu}).
Notice that $\de \Nthat_{af}$ contains
 the support of the cokernel of the morphism
$B_\lambda$ for any $\lambda\in \Lambda^+$, and equals this
support if
$\lambda\in \Lambda^{++}$ (see \ref{arr} for notation).

For a morphism $\phi:V\to L$ in a tensor category over a characteristic zero field and $d\in \Zet_{>0}$
one can form
the Koszul complex $0\to \Lambda^d(V) \to \Lambda^{d-1} (V) \otimes L \to \cdots
\to \Lambda^i(V)\otimes Sym^{d-i}(L) \to \cdots \to Sym^d(L) \to 0$. In the examples below some exterior
power of $V$ vanishes, and we will let $d$ be the maximal integer such that $\Lambda^d(V)\ne 0$,
the resulting complex will be called the Koszul complex associated to $\phi$.

Let
  $\bbK_{\lambda}\in Kom(Coh_{fr}^{\LG\times \LT}(\Nthat))$
 denote the  Koszul complex associated to  $B_\lambda$.
Thus
$$\bbK_\lambda=\left( 0\to \O= \Lambda^d(V_{\lambda}) \otimes \
\O \to \Lambda^{d-1}(V_{\lambda}) \otimes
\O (\lambda)
\to \cdots % \O((d-2)\lambda) \otimes \Lambda^2(V_\lambda)\to
\O((d-1)\lambda)\otimes V_\lambda \to \O(d\lambda)\to 0\right)$$
The key step is the following
\begin{Lem} We have $\Ftil(\bbK_\lambda)=0$ for all $\lambda\in \Lambda^+$.
\end{Lem}
\proofpt We keep the notations of the previous section. Thus
$\Ftil(\bbK_\lambda)$ is a complex of objects of $\A$. To see that
$\Ftil(\bbK_\lambda)$ is acyclic it is enough to see that $gr(\Ftil(\bbK_\lambda))\in Kom(\grA)$
is acyclic. The latter is a complex in $\grA\cong Rep(\LT)$. Since
the differential in $\bbK_\lambda$ is obtained from the arrow
$B_\lambda$ by tensoring with $V_\lambda$ and taking direct
summands, Theorem \ref{grfib} together with
Proposition \ref{moncen}(a) show that $gr(\bbK)\in
Kom(\grA)\cong Kom(Rep(\LT))$ is identified with the Koszul
complex associated to the non-zero map $ V_\lambda|_{\LT}\to \la$ in $Rep(\LT)$.
Since the latter complex is  acyclic, we get the statement. \epf

Now the definition of $F$ follows from the next

\begin{Lem} Let
$Hot_0(Coh_{fr}^{\LG\times\LT}(\Nthat_{af}))
\subset Hot(Coh_{fr}^{\LG\times \LT}(\Nthat_{af}))$ be the thick
subcategory of complexes whose cohomology is
supported on $\de \Nthat$.

a) Any $\F\in  Hot_0(Coh_{fr}^{\LG\times\LT}(\Nthat_{af}) )
\subset Hot(Coh_{fr}^{\LG\times \LT}(\Nthat_{af}))$ is a direct
summand in
$\F\otimes \bbK_\lambda$ for some $\lambda$.

b) The functor of restriction to $\Nthat$
provides an equivalence 
\begin{equation}\label{quot}
Hot(Coh_{fr}^{\LG\times\LT}(\Nthat_{af}))/
Hot_0(Coh_{fr}^{\LG\times\LT}(\Nthat_{af}))\iso D^{\LG\times \LT}(\Nthat)\cong
D^\LG(\Nt).
\end{equation}
\end{Lem}
\proofpt a)  It is clear that $Hot(Coh_{fr}^{\LG\times \LT}(\Nthat_{af}))$ is
identified with a full subcategory in
$D^b(Coh^{\LG\times\LT}(\Nthat_{af})$. For any  $\F\in
Hot_0(Coh_{fr}^{\LG\times\LT}(\Nthat_{af}))$
the corresponding object of $D^b(Coh^{\LG\times\LT}(\Nthat_{af})$
can be represented by
a finite complex of coherent sheaves set-theoretically supported on
$\de\Nthat$. This is clear by the following well-known fact (cf. e.g. \cite{izvrat}, Lemma 3(b)):
for an algebraic variety $X$ and a closed subvariety $Z\subset X$ the tautological
functor provides an equivalence  $D^b(Coh_Z(X))
\cong D^b_Z(Coh(X))$, where $Coh_Z(X)\subset Coh(X)$ is the full subcategory of sheaves set-theoretically supported on $Z$, and
 $D^b_Z(Coh(X))\subset D^b(Coh(X))$ is the full subcategory
of complexes whose cohomology sheaves lie in $Coh_Z(X)$. 

If $C$ is a finite 
complex of coherent sheaves set-theoretically supported on
$\de\Nthat$, then it is scheme-theoretically supported on some nilpotent neighborhood of 
$\de\Nthat$. 
For some $\lambda\in \Lambda^+$ the restriction of $B_\lambda$ to this neighborhood 
vanishes, thus
we have
$B_\lambda\otimes id_C=0$. Hence
$$\bbK_\lambda \otimes C\cong \oplusl_i \F\otimes
  \Lambda^i(V_{\lambda})\otimes
\O((d-i)\lambda),$$ which implies (a).

b) It suffices to check that the image of the functor \eqref{quot}
generates $D^{\LG}(\Nt)$ as a triangulated category; and that is a full
embedding. Here the first statement follows from Lemma \ref{Ogen}(a)
below; and the second one 
is a particular case of the following general statement. \epf

\begin{Sublem}\label{identify}Let $A$ be a finitely generated commutative
algebra graded by $\Zet_{\geq 0}^N$, and let $X$ be the
 corresponding multi-Proj scheme. Let $D_A^{\fr}$ be
the homotopy category of free $\Zet^N$-graded $A$-modules, and
$D_A^{\fr, 0}$
 be the full subcategory of complexes whose localization to
$D^b(Coh_X)$ is zero. Then $D_A^{\fr}/D_A^{\fr,0}$ is identified
with a full subcategory in $D^b(Coh_X)$.

Same is true for the categories of $H$-equivariant
sheaves/modules, where $H$ is a reductive algebraic group acting
on $A$.
\end{Sublem}

\proofpt For any finite complex $C\in D_A^{\fr}$, and any
$\lambda_0\in \Zet^N$ there exists $C'\in D_A^{\fr}$, and a
morphism $f:C'\to C$, such that $cone(f)\in D_A^{\fr,0}$, and $(C')^i$
is a sum of modules of the form $A\otimes V(\lambda)$,
$\lambda_0-\lambda\in \Zet_+^N$. (To see this pick  a
($H$-invariant) subspace $V\in A_\mu$, for $\mu$ large,
 such that $A/V\cdot A$ is
 supported on the complement to the cone over $X$, then consider
the Koszul complex
$$K=\left(0\to A(-d\mu)\to A(-(d-1)\mu)\otimes V^* \to
\cdots \to \Lambda^d(V^*)\otimes A \to 0 \right)$$
placed in degrees from $-d$ to 0.
 We have $K\otimes
C\in D_A^{\fr,0}$, and the complex
$Ker(K\otimes
\Lambda^d(V) \otimes C \to C)[-1]$ has the required form).

 Now given $B\in D_A^{\fr}$ we can find $\lambda_0$ such that
$Hom_{D_A}(A(\lambda)\otimes V,B)\iso
Hom_{D^b(Coh_X)}(\O(\lambda)\otimes V, \L(B))$ whenever
$\lambda_0-\lambda\in \Zet_+^n$ where $\L:D_A^{\fr}\to D^b(Coh_X)$
is the localization functor. Then also $Hom_{D_A^{\fr}}(C',B) \iso
Hom(\L(C'), \L(B))$ for $C'$ as above, which
%Since $$Hom_{
%D_A^{fr}/D_A^{fr,0}}(C,B)=\varinjlim Hom_{D_A^{fr}}(C',B)$$ where
%the limit is over the category of arrows $C'\to C$ whose cone is
%in $D_A^{fr,0}$, this
 implies the statement. \epf
\end{section}

\begin{section}{Proof of Theorem \ref{Phi}}
\subsection{Intermediate statements}

%The proof will be given at the end of this section.
We will deduce the Theorem from the next two statements.

Recall that $\kappa$ denotes the bijection $\Lambda\to \fW$. For
$\F\in \DIW$, and $\mu\in \Lambda$ set
$\Stalk_\mu(\F)=i_x^*(\F)[-\dim \Fl^{\kappa(\mu)}]$,
$Costalk_\mu(\F)=i_x^!(\F)[\dim \Fl^{\kappa(\mu)}]$ for $x\in
\Fl^{\kappa(\mu)}$; these are objects of the derived category of
vector spaces defined up to an isomorphism.

\begin{Prop}\label{tpW}
 For $k$ algebraically closed we have
$$\Stalk_\mu(\FIW(V))\cong \Ql^{\oplus [\mu :V]}\cong
\Costalk_\mu(\FIW(V)),$$ where $[\mu :V]$ is the multiplicity of
the weight $\mu$ in $V$.
\end{Prop}

\begin{Prop}\label{inje} For $V\in Rep(\LG)$, $\mu \in \Lambda^+$
the  map \begin{equation}\label{map}
Hom_{D^\LG(\Nt)}(V\otimes \O,
\O(\mu))\to Hom(\FIW(V\otimes \O), \FIW(\O(\mu)))
\end{equation}
is injective.
\end{Prop}

\subsection{Proof of Theorem \ref{descrIWiprec}} We now deduce
the Theorem from Propositions \ref{tpW}, \ref{inje}.

\begin{Lem}\label{Ogen} a) The objects $\O(\lambda)$, $\lambda\in
\Lambda$ generate $D^\LG(\Nt)$ as a triangulated category.

b) The objects $\O(-\lambda)\otimes V_\mu$, $\lambda,\mu\in
\Lambda^+$ generate $D^\LG(\Nt)$ as a triangulated category.
\end{Lem}

\proofpt (a) Since $\Nt$ is smooth every object of $D^\LG(\Nt)$ is
represented by a finite complex of $\LG$-equivariant vector
bundles. We now claim that every such vector bundle is filtered by
line bundles $\O(\lambda)$. Let $\E$ be such a vector bundle. It
is enough to show that there exists a $\LG$-equivariant injection
of vector bundles $\O(\lambda)\imbed \E$.

 We have an equivalence
$Coh^\LG(\Nt)\cong Coh^\LB(\Ln)$, $\E\mapsto \E|_{\Ln}$. Let
$M=\Gamma(\E|_{\Ln})$; then the data of an injection
$\O(\lambda)\imbed \E$ is equivalent to the data of an element
$v\in M$ such that $\LB$ acts on $v$ by the character $\lambda$,
and $v$ projects to a nonzero element in the coinvariants
$M/(\Ln)^*M$. It is easy to see that if $\lambda$ is a lowest
weight of $\LT\subset \LB$ in $M$ (which necessarily exists,
because the set of weights of $M$ is readily seen to be bounded
below) then every $v$ of weight $\lambda$ satisfies these
requirements.

(b) In view of statement (a) it is enough to show that for all
$\lambda$ the line bundle $\O(\lambda)$ lies in the triangulated
category generated by $\O(-\eta)\otimes V_\mu$, $\eta,\mu\in
\Lambda^+$. Acyclicity of the Koszul complex $\bbK_\mu$ (see
section \ref{fa}) shows that for all $\mu\in \Lambda^+$ the
 sheaf $\O(d\mu)$ lies in the triangulated category generated
by $\O((d-k)\mu)\otimes V$, $k\geq 1$, $V\in Rep(\LG)$; twisting
it by $\O(\lambda-d\mu)$ we see that $\O(\lambda)$ lies in the
triangulated category generated by $\O(\lambda-k\mu)\otimes V$,
$k\geq 1$, $V\in Rep(\LG)$. But given $\lambda \in \Lambda$ we can
find $\mu \in \Lambda^+$, such that $k\mu-\lambda\in  \Lambda^+$
for all $k\geq 1$. \epf

\begin{Lem}\label{Pgen}
The objects $\FIW(\O(\lambda))$, $\lambda\in \Lambda$
%(respectively, $\FIW^{gr}(\O(\lambda)(n))$, $\lambda\in \Lambda$,
%$n\in \Zet$)
 generate $\DIW$ % (respectively, $\DIW^{gr}$) 
as a triangulated category.
\end{Lem}

\proofpt Theorem \ref{IW}(b) implies that the support of $\FIW(\O(\lambda))$
is contained in the closure of $\Fl^{\kappa(\lambda)}$.

Furthermore, it is shown in \cite[4.1.2, Lemma 11]{cohotil} that the restriction of $J_\la$ to 
the $\GO$-orbit of $\Fl_\la$ coincides with the restriction of the standard sheaf $j_{\la!}$.
Thus Lemma \ref{nunu}(c) shows that the restriction of  $\FIW(\O(\lambda))=Av_\Psi(J_\la)$
to $\Fl^{\kappa(\lambda)}$ has rank 1. 
 This implies the Lemma. \epf

\subsubsection{Proof of Theorem \ref{descrIWiprec}}
We first check that $\FIW$ is a full embedding, i.e. that the map
\begin{equation}\label{full}
Hom^\bu _{D^\LG(\Nt)}(\F, \G)\to Hom^\bu _{D^\LG(\Nt)}(\FIW(\F),
\FIW(\G))
\end{equation}
is an isomorphism.

 It is known e.g. by results of \cite{KLT} that
for $\lambda\in \Lambda^+$ and $V\in Rep(\LG)$ we have
 $$Ext^{i}_{Coh^\LG(\Nt)}(V\otimes \O,
 \O(\lambda))=Hom_\LG(V,H^i(\Nt, \O(\lambda))=0$$ for
$i\ne 0$; and
$$\dim Hom_{Coh^\LG(\Nt)}(V\otimes \O,
 \O(\lambda))= [\lambda:V].$$
 The latter also equals $$\dim Hom(\FIW(V\otimes \O), \FIW(\O(\lambda))
= \dim Stalk_\lambda(\FIW(V\otimes \O))$$ by Proposition
\ref{tpW}. Thus Propositions \ref{tpW}, \ref{inje} imply that
\eqref{full} is an  isomorphism for $\F=V\otimes \O$,
$\G=\O(\lambda)$, $\lambda\in \Lambda^+$. Since $\FIW(\F\otimes
\O(\lambda))=\FIW(\F) *J_\lambda$, and the functor $\F \mapsto
\F\otimes \O(\lambda)$ is invertible we see that it is also
an isomorphism for $\F= V\otimes \O(-\lambda)$, $\G=\O$.
 Hence by Lemma \ref{Ogen}(b) it is an isomorphism for $\G=\O$ and
 all $\F$. Again twisting by $\O(\lambda)$ we deduce that it is an
 isomorphism for all $\F$ and $\G=\O(\lambda)$. Hence this is also
 true for all $\F$, $\G$ by Lemma \ref{Ogen}(a).

We proved that $\FIW$ is a full embedding.
It is then essentially
surjective by Lemma \ref{Pgen}. \epf

\subsection{Proof of Proposition \ref{inje}: category $\Po$ and the regular orbit}
Let  $\N^0\subset \N$ be the open $\LG$ orbit, and $N_0\in \N^0$
be an element. We introduce yet another auxiliary
category $\P_I^0$;  this category can be identified with
$Coh^\LG(\N^0)=Rep(Z_\LG(N_0))$.

 Let $D_I^{\ne 0}\subset D_I$ be the thick subcategory generated
 by $L_w$, $\ell(w)\ne 0$; and let $D_I^0$ be the quotient
 category. Let also $\PoI\subset D_I^0$ be the image of $\P_I$.
 Thus $\PoI=\P_I/(\P_I\cap D_I^{\ne 0})$ is an abelian category;
 %%Why??!!
  it has one irreducible object $L_e$ if
$G$ is simply connected; in general the number of irreducible
objects in $\Po_I$ equals $\#\pi_1(G)$.

The convolution of any $X\in D_I^{\ne 0}$ with any object of
$D_I$ lies in $D_I^{\ne 0}$; thus the convolution induces a
monoidal structure on $D_I^0$ (which we denote by the same
symbol). Moreover, the abelian subcategory $\PoI\subset D_I^0$ is
monoidal.

Let $F_0:Rep(\LG)\to \PoI$ be the composition of the projection
functor $\P_I\to \PoI$ with  $\Z$. Then $F_0$ respects the
monoidal structure. Let $\cT$ denote the full subcategory of $\PoI$
consisting of all
 subquotients of $F_0(V)$ for some $V\in Rep(\LG)$ (in fact, one can show
  that $\cT=\PoI$, but we neither use nor prove this here).
  Also $\M$ induces a tensor endomorphism of the functor $F_0$.

Let us recall that for any subgroup $H\subset Z_\LG(N_0)$ the
restriction functor $res^{\LG}_H$ carries a canonical tensor
endomorphism induced by $N_0$.

\begin{Lem}\label{fibA}
  There exists a subgroup $H\subset Z_\LG(N_0)$, and an
  equivalence of monoidal categories $\cT\cong Rep(H)$, which
  intertwines $F_0$ with the restriction functor $res^{\LG}_H$,
  and sends the tensor endomorphism $\M$ into the endomorphism
  induced by $N_0$.
\end{Lem}

\proofpt See \cite{B}. \epf

\begin{Rem}  It is not difficult to deduce from the results
of the present paper that in fact $H=Z_\LG(N_0)$;
and  also  that $\cT=\PI^0$.
\end{Rem}

\begin{Rem} Lemma \ref{fibA} is a particular case of a more
general result proved in \cite{B}, which relates representations
of a centralizer of any nilpotent in $\Lg$ to a two-sided cell in
$W$. This result is deduced from some "non-elementary" Theorems of
Lusztig about (asymptotic) Hecke algebras, see \cite{cells4}
(which rely on the theory of character sheaves).  However, the
particular case used in Lemma \ref{fibA} depends only on the
elementary particular case of Lusztig's Theorems,
which deal with the maximal cell and the regular
nilpotent orbit (the corresponding fact about the Hecke algebra amounts
to the computation of the action of the center of the affine Hecke
algebra in the Steinberg representation).
\end{Rem}

For $V_1,V_2\in Rep(\LG)$ Lemma \ref{fibA} yields an injective map
\begin{equation}\label{mapone}
Hom_{Z_\LG(N_0)}(V_1,V_2)\to Hom_{\PoI}(F_0(V_1),F_0(V_2)).
\end{equation}
On the other hand, we have
$$Hom_{Z_\LG(N_0)}(V_1,V_2)=Hom_{Coh^\LG(\N^0)}(V_1\otimes \O,
V_2\otimes \O)=Hom_{Coh^\LG(\Nt)}(V_1\otimes \O, V_2\otimes \O),$$
thus the functor $F$ induces another map
\begin{equation}\label{maptwo}Hom_{Z_\LG(N_0)}(V_1,V_2)\to Hom_{\P_I}(\Z(V_1),\Z(V_2))
\to Hom _{\P_I^0}(F_0(V_1),F_0(V_2)).
\end{equation}

\begin{Lem}
The map \eqref{mapone} coincides with the composition in \eqref{maptwo}.
\end{Lem}

\proofpt This follows from the uniqueness statement in Proposition
\ref{FtoC}(a), since both \eqref{mapone} and \eqref{maptwo} send
the "tautological" tensor endomorphism to the logarithm of
monodromy endomorphism $\M$. \epf

\begin{Cor}\label{in}
The composed map \eqref{maptwo} is injective. \epf
\end{Cor}

\subsubsection{Proof of Proposition \ref{inje}}
We first claim that for any $\lambda\in \Lambda^+$ the sheaf
$\O(\lambda)\in Coh^\LG(\Nt)$ can be realized as a subsheaf in
$V\otimes \O_\Nt$ for some $V\in Rep(\LG)$. Indeed, for a simple
coroot (root of $\LG$) $\alpha\in \Lambda$ let us denote by
$D_\alpha \subset \Nt$ the $\LG$-invariant divisor
$T^*(\LG/\LP_\alpha)\times _{\LG/\LP_\alpha} \LG/\LB$, where
$\LP_\alpha\subset \LG$ is the corresponding minimal parabolic.
Then it is easy to see that
$$\O_{\Nt}(-\alpha)\cong \O(D_\alpha);
$$
thus we have an injective map of sheaves $\O\imbed \O(-\alpha)$.
Taking tensor products we get also injections $\O(\lambda)\imbed
\O(\lambda- 2n \rho)$. For large $n$ we have $\lambda-n\rho \in
-\Lambda^+$, so we get an injection $$\O(\lambda-n\rho)\imbed
V_{w_0(\lambda-n\rho)}\otimes \O.$$

Thus it suffices to see that the map
\begin{equation}\label{mapthree}
Hom(V_1\otimes \O, V_2\otimes \O)\to Hom (\FIW(V_1\otimes \O),
\FIW(V_2\otimes \O))
\end{equation}
is into. The functor
$Av_\Psi:\P_I\to \PIW$ is exact by
 Proposition \ref{IWisfPprop}(a); by Lemma \ref{nunu}(a)
it does not kill $L_w$ for $w\in \fW$, in particular
for $\ell(w)=0$. Hence it does not kill any
morphism  whose image in $\PoI$ is nonzero. The statement now
follows from Corollary \ref{in}. \epf

\subsection{Proof of Proposition \ref{tpW}} We first prove the
following

\begin{Thm}\label{tilt} $\Stalk_\mu(\FIW(V))$, $\Costalk_\mu(\FIW(V))$
are concentrated in homological degree 0 for all $\mu\in
\Lambda^+$, $V\in Rep(\LG)$.
\end{Thm}
\begin{Rem}\label{remtiltasph}
Theorem \ref{tilt} says that $\FIW(V)$ is a {\it tilting} object
of $\PIW$.
\end{Rem}
\begin{Rem}
We do not know whether the following strengthening of Theorem
\ref{tilt} is true: ``for every  convolution exact object
$\F$ of $\P_I$ the sheaf $\Delta_0*\F$ is tilting'' (cf. Remark
\ref{remtilt} in section \ref{unifi}).
\end{Rem}

\begin{Rem}\label{Koszul}
Recall that the {\it parabolic-singular} Koszul duality
is an equivalence between the mixed versions of parabolic and
singular categories $O$, see \cite{BGS}. An appropriate version of this
equivalence (see \cite{BG}) sends irreducible objects into tilting ones.
 The parabolic category
$O$ is equivalent to the category of perverse sheaves on 
the partial flag variety $G/P$.
Using (a variation of) the result of \cite{MS} one can realize 
the singular category as an appropriate category
of Whittaker sheaves. One can try to generalize this picture by
replacing $G$ by the loop group $\GK$, and $P$ by 
the maximal parahoric $\GO$.
Thus we are led to the conjecture that there exists an equivalence
between the mixed versions of $\DIW$ and the category of Iwahori
monodromic sheaves on the affine Grassmanian. In fact, this conjecture
can be derived from a combination of the results of this paper and
those of \cite{ABG}, or by adapting the method of \cite{BGS}; see also discussion in \cite[1.2]{cohotil}.

In view of some
formal properties of this duality (in particular, the fact
that the central sheaves are  Koszul self-dual\footnote{This fact
was suggested to us by M.~Finkelberg.}) the statement of
Theorem \ref{tilt} is Koszul dual  to the statement that
 $\pi_*(Z_\lambda^{Weil})$ is  simple of weight 0; the
latter statement is clear from the definition
of $\Z$ together with the fact that nearby cycles commute
with proper direct image, cf. \cite{KGB}, Theorem 1(d).
\end{Rem}

The Theorem will be deduced from the following two statements.

\begin{Lem}\label{ifthen} If Theorem \ref{tilt} holds for two
representations $V_1,V_2$, then it holds for $V=V_1\otimes V_2$.
\end{Lem}

 \begin{Lem}\label{miqmi}  Theorem \ref{tilt} holds if $V=V_\lambda$, where
 $\lambda \in \Lambda^+$ is either minuscule or  quasi-minuscule 
(i.e. is the short dominant root).
 \end{Lem}

\subsubsection{Proof of Theorem \ref{tilt}} If $\LG$ is not adjoint let $V$
be the sum of its minuscule irreducible representations; otherwise
let $V$ be the quasi-minuscule representation. Lemma \ref{miqmi}
shows that the statement of the Theorem holds for $V$. However,
it is easy to see $V$ is a faithful representation; hence it induces
a surjective map from functions on $End(V)$ to functions on $\LG$.
Since $V$ is self-dual, 
any irreducible 
representation of $\LG$ is a direct summand of $V^{\otimes n}$
for some $n$. Thus the Theorem follows by
Lemma \ref{ifthen}. \epf

\begin{Rem}
The trick of reduction to the special case of a (quasi-)minuscule
representation was also (independently) used in \cite{PN}.
\end{Rem}

\subsubsection{Proof of Proposition \ref{tpW}}
The Proposition follows from Theorem \ref{tilt} and the following
Lemma.

\begin{Lem}\label{eul}
The Euler characteristic $\sum (-1)^i \dim H^i(Stalk_\mu
\FIW(V\otimes \O))$ equals $[\mu:V]$, the multiplicity of the
weight $\mu$ in $V$.
\end{Lem}

\proofpt For  a triangulated category $D$ we will denote its
Grothendieck group by $K(D)$, and for $X\in D$ will let  $[X]\in
K(D)$ be its class. We have an isomorphism $K(D_I)\cong \Zet[W]$,
$[j_{w!}]\mapsto w$. This isomorphism is compatible with the
algebra structure, where the one on $K(D_I)$ comes from the
convolution on $D_I$; this follows from the equalities
$$\begin{array}{ll}
&[j_{w_1!}]\cdot  [j_{w_2!}]=[j_{w_1!} * j_{w_2!}]=[j_{w_1w_2!}] \
\ {\rm for} \ \
\ell(w_1w_2)=\ell(w_1)+\ell(w_2);\\
&[j_{s!}]^{-1}=[j_{s*}]=[j_{s!}],
\end{array}$$
where $s$ is a simple reflection. It follows that $[J_\lambda]=
[j_{\lambda!}]$ for all $\lambda \in \Lambda$. Thus Theorem
\ref{filtrWakr} implies that $$[Z(V)]=\oplusl_\mu [\mu:V]\cdot
[J_\mu]=\oplusl_\mu [\mu:V]\cdot [j_{\mu!}].$$ Thus $$
[\FIW(V\otimes \O)]=\oplusl_\mu [\mu:V]\cdot
[\Delta_0*j_{\mu!}]=\oplusl_\mu [\mu:V]\cdot [\Delta_\mu],$$
which implies the statement of the Lemma. \epf

\subsubsection{Proof of Lemma \ref{ifthen}}  We will consider the  condition on stalks; the
one on costalks is treated similarly. Notice that this condition
is equivalent to saying that $\FIW(V)$ carries a filtration with
subquotients isomorphic to  $\Delta_\mu$. If this is the case for
$V=V_1$, then $\FIW(V_1\otimes V_2)$ carries a filtration with
subquotients of the form
\begin{equation}\label{subq}
\Delta_\mu*\Z(V_2)\cong \Delta_0*\Z(V_2)*j_{\mu!}, \end{equation}
where the central property of the sheaf $\Z(V_2)$ is used. The
Lemma will be proven if we show the corresponding statement for
the stalk of the sheaf $\F$ appearing in either side of
\eqref{subq}.  $\F$ is a perverse sheaf by exactness of
convolution with $\Z(V_2)$; hence $\Stalk_\mu(\F)\in D^{\leq 0}$
by the definition of a perverse sheaf. The opposite estimate
follows from the assumption that $\Delta_0*\Z(V_2)$ has a
filtration with subquotients $\Delta_\nu$, and the following \epf

\begin{Sublem}
For all $\lambda,\nu\in \Lambda$, $w\in W$ we have
$$\Stalk_\lambda(\Delta_\nu*j_{w!})\in D^{\geq 0}.$$
\end{Sublem}
\proofpt The statement is equivalent to
%\begin{equation}\label{svone}
$$
\Delta_\nu*j_{w!}\in \langle \Delta_\mu[i]\ ,\ \ i\leq 0, \mu\in
\Lambda \rangle
$$
%\end{equation}
(cf. sections \ref{prooW}, \ref{genizvr} for notations). Since
$\Delta_\lambda \cong \Delta_0*j_{u!}$ for any $u\in W_f\cdot
\lambda$ by Lemma \ref{nunu}(c), the latter follows from the following
statement (see e.g. \cite{intertwi})
\begin{equation}\label{svtwo}
j_{w_1!}*j_{w_2!}\in \langle j_{w!}[i]\ ,\ \ i\leq 0, w\in
W\rangle.
\end{equation}
To verify \eqref{svtwo} we can assume that $w_2=s$ is a simple
reflection. If $\ell(w_1\cdot s)>\ell(w_1)$, then
$j_{w_1!}*j_{s!}\cong j_{w_1\cdot s!}$, so \eqref{svtwo} is clear.
If $\ell(w_1\cdot s)< \ell(w_1)$, then we have an exact triangle
$$j_{w_1!}\oplus j_{w_1!}[-1]\to j_{w_1!}*j_{s!}\to j_{w_1s!},$$
which shows \eqref{svtwo} in this case also. \epf

\subsubsection{Proof of Lemma \ref{miqmi}}
\begin{Lem}\label{indepe}
 For $w_f\in W_f$, $V\in Rep(\LG)$ we have an isomorphism
$$\Stalk_\lambda(\FIW(V\otimes \O))\cong \Stalk_{w_f(\lambda)}(\FIW(V\otimes \O)).$$
\end{Lem}

\proofpt Let $s\in W_f$ be a simple reflection. Then
$\Delta_0*L_s=0$ by Lemma \ref{nunu}(a) above. By the central
property of $\Z(V)$ we have also
$$\FIW(V\otimes \O)*L_s=\Delta_0*L_s*\Z(V)=0.$$
For $X\in \DIW$, and $\lambda\in \Lambda$ such that
$s(\lambda)\preceq \lambda$ it is easy to construct an exact
triangle $$\Stalk_{\lambda}(X)[-1]\to \Stalk_\lambda (X*L_s) \to
\Stalk_{s(\lambda)}(X)\to \Stalk_\lambda(X).$$ Thus
$\Stalk_\lambda(X)\cong \Stalk_{s(\lambda)}(X)$ provided that
$X*L_s=0$. This proves the statement of the Lemma for $w_f=s$, and
hence for all $w_f\in W_f$. \epf

\begin{Lem}\label{dV} For $V\in Rep(\LG)$ let $d_V=\dim V^{N_0}$,
  where $N_0$ is a
regular nilpotent element. Then we have
\begin{equation}\label{eqocen}\dim Hom(\Delta_0*\Z(V),
\Delta_0)\leq d_V.
\end{equation}
\end{Lem}

The proof of the Lemma will rely on the following result of
D.~Gaitsgory (unpublished); we reproduce the proof in the
Appendix.
\begin{Thm}\label{monvra}
There exists an element $\mon={\mon_\F}$ of the center of $\P_I$
such that

i) $\mon_{\Z(V)}=\M_{V}$ for $V\in Rep(\LG)$.

ii) $\mon_L=0$ for any irreducible object $L\in \P_I$.
\end{Thm}

\proof of Lemma \ref{dV}. We have $$Hom(\Delta_0*\Z(V), \Delta_0)
=Hom_{\fP_I}(\Z(V), L_e)=Hom_{\fP_I} (\Z(V)_{\M_V},L_e)\imbed
Hom_{\PoI}(\Z(V)_{\M_V},L_e),$$ where the first equality follows
from  Proposition \ref{IWisfPprop}(b); and the second one from
Theorem \ref{monvra}. But Lemma \ref{fibA} implies that
$\Z(V)_{\M_V}\mod \PI^{\ne 0}$ has length $d_V$, which shows that
$\dim Hom_{\PoI}(\Z(V)_{\M_V},L_e)\leq d_V$. \epf

{\it Proof of Lemma \ref{miqmi}} (conclusion). It follows from
Lemma \ref{bFl} that $\Fl^{\lambda}$ is open in the support
of $\FIW(V_\lambda)$, and
%%eto mini-maksnaya kletka, poetomu \pi ot nee imeet bol'shuyu r-t'
$$\Stalk_{w_0(\lambda)}(\FIW(V_\lambda\otimes \O))\cong \Ql.$$
By Lemma \ref{indepe} we conclude that
\begin{equation}\label{onU}
Stalk_{w_f(\lambda)}(\FIW(V_\lambda\otimes \O))\cong \Ql\ \ \ \
{\rm for} \ \ w_f\in W_f.
\end{equation}
 If
$\lambda$ is minuscule then \eqref{onU} implies the statement of
the Proposition, because in this case the support of
$\FIW(V_\lambda\otimes \O)$ contains $\Fl^\mu$ iff $\mu \in
W_f(\lambda)$.

 Assume now that $\lambda$ is
quasi-minuscule, i.e. $\lambda$ is the short dominant coroot. Then
the support of $\FIW(V_\lambda\otimes \O)$ contains $\Fl^\mu$ iff
either $\mu \in W_f(\lambda)$ or $\mu=0$. The first case is
treated by \eqref{onU}, so it remains to consider the case
$\mu=0$.

We first  claim that $Stalk_0(\Delta_0*\Z(V))$ can only be
concentrated in degrees 0 and $-1$. This follows from the exact
triangle
$$j_!j^*( \Delta_0*\Z(V))\to  \Delta_0*\Z(V)\to
i_*i^*(\Delta_0*\Z(V)),$$ where $i$ is the embedding of $G/B\imbed
\Fl$, and $j$ is the embedding of its complement. \eqref{onU}
shows that $j^*( \Delta_0*\Z(V))$ carries a filtration whose
subquotients are of the form $j^*(\Delta_\mu)$. Since
$j_!j^*(\Delta_\mu)=\Delta_\mu$ for $\mu \ne 0$  we see that
$j_!j^*( \Delta_0*\Z(V))$ is a perverse sheaf, hence
$i_*i^*(\Delta_0*\Z(V))$ is concentrated in homological degrees
$0$ and $-1$.

To finish the proof it now suffices to check that
\begin{equation}\label{nrvo}
 \dim H^0(Stalk_0(\Delta_0*\Z(V)))\leq \sum (-1)^i \dim
H^i(Stalk_0(\Delta_0*\Z(V))).
\end{equation}
 Here the vector
space in the left hand side is dual to the vector space
$Hom(\FIW(V\otimes \O), \Delta_0)$, so its dimension is estimated
 by Lemma \ref{dV}. The right hand side of
\eqref{nrvo} is computed in Lemma \ref{eul}. To see \eqref{nrvo}
it remains to notice that for $\lambda$ quasi-minuscule we have
$$[0:V_\lambda]=\dim V_\lambda^h=\dim V_\lambda^{N_0},$$
where $h$ is the semisimple element in a regular $sl(2)$ triple
$e=(N_0, h,f)$; here the first (respectively, the second)
equality is true 
 because every non-zero weight of $V_\lambda$ is
 a
root, hence is non-zero (respectively, even)
on the Cartan of a principle $SL(2)$.
 \epf

\subsection{Proof of Theorem \ref{IWisfP}}\label{IWisfPkonec}
%We describe the unmixed case.
In view of Proposition
\ref{IWisfPprop} it suffices to check that $Av_\Psi:\fP_I\to \PIW$
% $Av_\Psi^{gr}:\fP_I^{gr}\to \PIW^{gr}$
is essentially surjective. By Theorem \ref{descrIWiprec} any $X\in
\PIW$ is isomorphic to $\FIW(\F)=pr_f(F(\F))$ for some $\F\in
D(\Nt)$. Since $Av_\Psi$ is exact by Proposition \ref{IWisfPprop}(a) we
see that $X\cong Av_\Psi(F(\F))\cong Av_\Psi(H^{0,p}(F(\F))$,
which shows essential surjectivity of $Av_\Psi$. 
%The argument for $Av_\Psi^{gr}$ is parallel. 
\epf

\subsection{Application: Whittaker integrals of central sheaves}
For $\F\in D$ set
 $Wh_w^i(\F)= H_c^i(\Psi_w^*(\AS)\otimes
\iota_w^*(\F))$. If $\F$ is endowed
with a Weil structure then  $Wh_w^i(\F)$ carries an action of Frobenius.

\begin{Thm}\label{Whitint}
 Let $\lambda\in \Lambda^+$, $w=\kappa(\mu)\in ^fW$.
 Recall that $ [\mu:V_\lambda]$ denote the multiplicity of the
weight $\mu$ in $V_\lambda$. Then we have

$$Wh_w^i(Z_\lambda)=0\ \ \ {\rm for}\ \ \ i\ne \dim(\Fl^w);$$
 $$\dim Wh_w^{\ell(w)}(Z_\lambda)= [\mu:V_\lambda].$$

Moreover, we have
$$Tr(Fr,Wh_w^{\dim(\Fl^w)}(Z_\lambda^{Weil}))=Q_{\lambda, \mu}(q^{1/2});$$
here the polynomial $Q_{\lambda,\mu}(t)$ is defined by 
\begin{equation}\label{ra}
Q_{\lambda,\mu}=Q_{\lambda,w(\mu)} \ \ \ \ \ {\rm for}\ \ \ \ w\in W_f;
\end{equation}
\begin{equation}\label{QP}
Q_{\lambda,\mu}=t^{\ell(\lambda)+\ell(w_0)}P_{\lambda,\mu}(t^2),
\end{equation}
where
$P_{\lambda,\mu}$ is the Kazhdan-Lusztig polynomial 
(the $q$-analogue of weight multiplicity), see
\cite{doistor}.
\end{Thm}

We will abuse notations by writing $\nabla_w$ for the Weil
sheaf $i_{w*}\psi_w^*(\AS)[\ell(w)](\lT)$. 

\begin{Lem}\label{avpsi}
a) For $\F\in D_I$ we have a canonical isomorphism
\begin{equation}\label{Whw}
Wh_w^{i}(Av_\Psi(\F))\cong Wh_w^{i+\ell(w_0)}
\left( \F(\frac{\ell(w_0)}{2})
\right).
\end{equation}
If $\F$ is equipped with a Weil structure, the isomorphism
is compatible with the Frobenius action.

b) For $\F\in \DIW$ we have \begin{equation}\label{Whwtwo}
Wh_w^{i}(\F)=Hom\left(\F,
\nabla_w[i+\ell(w)](\frac{\ell(w)}{2})\right);
\end{equation}
 for $\F$  equipped with the
Weil structure the isomorphism is compatible with the Frobenius action.
\end{Lem}

\proofpt It is easy to see that $\Fl^e$ is an orbit of the group
$\bI^-\cap \NFm$. It follows that $$Av_\Psi(\F)\cong Av_{\bI^-\cap
\NFm, \psi}^![\ell(w_0)](\frac{\ell(w_0)}{2}),$$ where
\begin{equation}\label{avpsieq}
 Av_{\bI^-\cap
\NFm, \psi}^!=a_!(\psi_{\bI^-\cap \NFm/U}^*(\AS)\otimes
pr^!(\F)), \end{equation}
 and $pr:\bI^-\cap \NFm/U
\times supp(\F) \to supp(\F)$, $a:\bI^-\cap \NFm/U \times
supp(\F)\to \Fl$ are the projection to the second factor, and the
action map respectively, while $U\subset \bI^-\cap \NFm/U $ is
some open subgroup which stabilizes all $x\in supp(\F)$. Proper
base change shows that
$$Wh_w^{i}(\F)=Wh_w^i(a_!(\psi_{H/U}^*(\AS)\otimes
pr^!(\F))$$ for any group subscheme $H\subset \NFm$, in evident
notations. This proves (a).

To prove (b)  notice that Lemma \ref{fW} implies that
$Wh_w^i(\Delta_{w'})=0$ unless $w= w'$, $i=\dim \Fl^w$; and
$$Wh_w^{\dim \Fl^w}(\Delta_{w})=\Ql(-\frac{\dim \Fl^w}{2}).$$
It follows that $$Wh_w^*(\F)=Wh_w^*(i_{w!}i_w^*(\F))= Hom^*(\F,
\nabla_w)[-\dim \Fl^w](-\frac{\dim \Fl^w}{2}).\ \ \ \epf$$

\proof of Theorem \ref{Whitint} (sketch).
The first two statements of the Theorem follow from
Lemma \ref{avpsi}(b), Theorem \ref{tilt} and Lemma \ref{eul}.

The equality \eqref{ra} follows from 
$$Wh_\mu^\bu(Z_\lambda*L_s)=0$$
which yields an isomorphism 
$$Wh_\mu(Z_\lambda)\cong Wh_{s(\mu)}(Z_\lambda)(-\frac{1}{2})$$
where $s\in W_f$ is a simple reflection such that $s(\mu)\prec \mu$.
 
 Let us prove \eqref{QP}.
Using Lemma \ref{avpsi} and
equivalence of Theorem \ref{descrIWiprec} we see that
$$Wh_\mu^{\dim \Fl}(Z_\lambda)\cong Hom
(V_\lambda\otimes \O,
\O(\mu)), $$
 %(n-\ell(\mu)-\ell(w_0)),$$
%= Hom (V_\lambda\otimes \O,\O(\mu)(n-\dim \Fl^w),$$
where the action of Frobenius on the left hand side corresponds
to $q^{(-\ell(\mu) -\ell(w_0))/2}$ times
 the action of the automorphism induces by the dilatation by $q$.

By the well-known interpretation
(due to Hesselink)
of $P_{\lambda,\mu}$ in terms of the character of the space of functions
on $\Nt$ (see e.g. \cite{Bry}, Lemma 6.1) the trace of the latter automorphism
equals
$ q^{(\ell(\lambda)+\ell(w_0))/2}P_{\lambda,\mu}(q)
$. \epf

 \end{section}

\begin{section}{Appendix}

\centerline{by ROMAN BEZRUKAVNIKOV and 
IVAN MIRKOVI\' C\footnote{Partly supported by an NSF grant.}}

\bigskip

\subsection{Sketch of proof of Theorem \ref{IW}}
It is easy to see that the morphisms
 $m_l^w:\Fl_w\underset{\bI}{\times}\Fl \to \Fl$
 $m_r^w:\Fl\underset{\bI}{\times}\Fl_w \to \Fl$ (restrictions of the
convolution map) are affine (notations
of section \ref{nota}). Hence for $w=(\lambda)^{-1}\cdot \mu \cdot w_f$,
$\lambda,\mu\in \Lambda^+$, $w_f\in W_f$ we have
$$J_w=j_{-\lambda!}*j_{\mu\cdot w_f*}=(m_l^{-\lambda})_!(
j_{-\lambda!}\underset{\bI}{\boxtimes}j_{\mu\cdot w_f*}
|_{\Fl_w\underset{\bI}{\times}\Fl})\in D^{p,\geq 0}$$
because $!$ direct image under an affine morphism is left exact in the
perverse $t$-structure (see \cite{BBD}, 4.1.2).
On the other hand
$$J_w=(m_r^{\mu\cdot w_f})_*(
j_{-\lambda!}\underset{\bI}{\boxtimes}j_{\mu\cdot w_f*}
|_{\Fl_w\underset{\bI}{\times}\Fl})\in D^{p,\leq 0}$$
because $*$ direct image under an affine morphism is right exact in the
perverse $t$-structure (see \cite{BBD}, 4.1.1).
The two observations together give  statement (a).

In view of (a), statement (b) will be proven if we check that the
Euler characteristic of the stalk $Eul_w(J_{w'})= \sum (-1)^i\dim
H^i(j_w^*(J_w'))$ is zero for $w'\not \preceq w$, and is nonzero
for $w'=w$. (Indeed, we have $Eul_w(J_{w'})\ne 0$ whenever $\Fl_w$
is open in the support of $J_{w'}$). In the proof of Lemma
\ref{eul} we have seen  the equality $[J_w]=[j_{w!}]$
 in the Grothendieck group $K(D_I)$.
It implies that $Eul_w(J_{w'})=(-1)^{\ell(w)}\delta_{w,w'}$. \epf

\subsection{Sketch of proof of Theorem \ref{monvra}}
Let $X$ be a scheme with a $\Gm$ action. Recall the notion of a
{\it monodromic} constructible complex on $X$, and the {\it
monodromy action} of the tame fundamental
group of $\Gm$ on such a sheaf, see \cite{V}. In fact, the definition of
{\it loc. cit.} works in the set up when $X$ is a cone over a base
scheme $S$, i.e. a closed $\Gm$-invariant subscheme in ${\mathbb
A}^n_S$. For an arbitrary $X$ we will say that a constructible
complex $\F$ on $X$ is monodromic if $j_!a^*(\F)$ is monodromic in
the sense of \cite{V}, where $a:\Gm\times X\to X$ is the action
map, and $j:\Gm\times X\imbed \Aone \times X$ is the embedding. If
$\F$ is a perverse sheaf, then we have $End(\F)\iso End
(j_!a^*(\F))$; so the monodromy action 
 on  $j_!a^*(\F)$ introduced in \cite{V}  yields
an action on $\F$,
 which we also call the $\Gm$-monodromy action.
If this action is unipotent, 
it  defines % (at least if it is unipotent)
 the {\it logarithm of monodromy} operator
$\F\to \F(-1)$, see e.g.
 \cite{Weil2}, \S 1.7.2, which we denote by $\mon_\F^{\Gm}$.

 Let us now describe the element $\mon$.
  Recall
that the pro-algebraic group $\AutO$ of automorphisms of $O$  acts
on $\Fl$. All $\bI$-orbits on $\Fl$ are $\AutO$ invariant. This
implies that $L_w$ are $\AutO$ equivariant; in particular, they
are equivariant with respect to the subgroup of dilations
$\Gm\subset \AutO$. Hence any $\F\in \P_I$ is monodromic with
respect to this action, and the monodromy action on $\F$ is
unipotent.

We set $\mon_\F=-\mon _\F^{\Gm}$.

 Then property (ii) is clear. To
establish (i) we recall that the functor $\Z(V)$ is defined as the
nearby cycles of a certain sheaf $\F_V$ on the space $\Fl_X$
defined in terms of a global curve $X$. We can assume that
$X=\Aone$, so that there is an  action of $\Gm\subset \AutO$ on
$X$ compatible with the action on $O$. It also induces an action
of $\Gm$ on $\Fl_X$. Then it is easy to see that the sheaf $\F_V$
is equivariant with respect to this action of $\Gm$. Now property
(i) follows from the following general

\begin{Claim} Let $X$ be a variety with a $\Gm$-action, and
$f:X\to \Aone$ be a function, such that $f(tx)=tf(x)$ for $x\in
X$, $t\in \Gm$. If $\F$ is a $\Gm$ equivariant perverse sheaf on
$X$, then the nearby cycles complex $\Psi_{f}(\F)$ is
$\Gm$-monodromic. Moreover, the monodromy action on $\Psi_f(\F)$
as on the nearby cycles sheaf factors through the tame fundamental
group of $\Gm$; the resulting action of this tame fundamental
group on $\F$ is opposite to the $\Gm$-monodromy action.
\end{Claim}

\proofpt
 This is a restatement of  \cite{V}, Proposition 7.1(a). More precisely,
the set up of {\it loc. cit.} is as follows. One considers a
constructible complex $\tilde\F=pr_1^*(\F)$ on $X\times \Gm$ where
$pr_1:\Gm\times X\to X$ is the projection; and a function $\tilde
f$ on $\Gm\times X$ given by $\tilde f(t,x)=tf(x)$. It is then proved
that the $\Gm$ monodromy action on $\Psi_{\tilde f}(pr_1^*(\F))
$ is opposite the canonical monodromy acting on the nearby cycles
sheaf. 
 Under our
assumptions $pr_1^*(\F)\cong a^*(\F)$, $\tilde f = a^*(f)$, where
$a:\Gm\times X\to X$ is the action map. Thus the statement of {\it
loc. cit.} implies the Claim. \epf

\begin{Rem}
This Claim implies that the action of monodromy on the functor
$\Z$ is unipotent (because the $\Gm$ monodromy is obviously unipotent,
as it is functorial, and is trivial on irreducible perverse sheaves).
This fact was proved in \cite{KGB} by a different argument.
\end{Rem}

\end{section}

\footnotesize{
{\bf S.A.}: $\,$  Department of Mathematics,
University of Toronto,
Toronto, Ontario, Canada M5S 2E4\\
%\hphantom{x}\ab\, 
{\tt hippie@math.toronto.edu}}

\footnotesize{
{\bf R.B.}: Department of Mathematics,
 Massachusetts Institute of Technology,
Cambridge MA,
02139, USA;\\ 
%\hphantom{x}\ab\, 
{\tt bezrukav@math.mit.edu}}

\footnotesize{
{\bf I.M.}: Department of Mathematics and Statistics,
University of Massachusetts,
 Amherst, MA 01003,
 USA;\\ 
%\hphantom{x}\ab\, 
{\tt mirkovic@math.umass.edu}}

\end{document}